\documentclass{amsart}
\usepackage[a4paper,hmargin=26mm,vmargin=28mm]{geometry}
\usepackage{amssymb}
\usepackage{etoolbox}
\usepackage{mathrsfs}
\usepackage{mathtools}
\usepackage{stmaryrd}
\usepackage{times}
\usepackage{xparse}

\usepackage{tikz}
\usetikzlibrary{matrix}

\usepackage{booktabs}
\usepackage[a4paper=true,pagebackref=true]{hyperref}
\renewcommand*{\backref}[1]{}
\renewcommand*{\backrefalt}[4]{%
  \ifcase #1 No citations.% shouldn't happen...
  \or [Page~#2.]%
  \else [Pages~#2.]%
  \fi%
}

\hypersetup{pdftitle={The irreducible characters of the alternating Hecke algebras},
  pdfauthor={Andrew Mathas and Leah Neves},
  colorlinks = true,
  linkcolor =blue,
  anchorcolor = red,
  citecolor = blue,
  urlcolor = blue
}
\usepackage{amsmath}

\numberwithin{equation}{section}
\swapnumbers

\def\NewTheorem#1{%
  \newtheorem{#1}[equation]{#1}
}
\usepackage[capitalise,nameinlink]{cleveref}
\crefname{equation}{}{}

\NewTheorem{Proposition}
\NewTheorem{Theorem}
\NewTheorem{Corollary}
\NewTheorem{Lemma}
\NewTheorem{Definition}
\theoremstyle{definition}
\NewTheorem{Example}
\AtEndEnvironment{Example}{\null\hfill$\Diamond$}%
\theoremstyle{remark}
\NewTheorem{Remark}

\newcount\tableauRow\newcount\tableauCol
\newcommand\Tableau[2][-12]{
  \begin{tikzpicture}[scale=0.38,draw/.append style={thick,black},baseline=#1mm]
    \tableauRow=0
    \foreach \Row in {#2} {
       \tableauCol=1
       \foreach\k in \Row {
          \draw(\the\tableauCol,\the\tableauRow)+(-.5,-.5)rectangle++(.5,.5);
          \draw(\the\tableauCol,\the\tableauRow)node{\small$\k$};
          \global\advance\tableauCol by 1
       }
       \global\advance\tableauRow by -1
    }
  \end{tikzpicture}
}

\newdimen\shadedBaseline\shadedBaseline=-4mm
\newcommand\ShadedTableau[2][]{
  \begin{tikzpicture}[scale=0.5,draw/.append style={thick,black},baseline=-4mm]
    \foreach\bx in {#1} { \filldraw[blue!20]\bx+(-.5,-.5)rectangle++(.5,.5); }
    \tableauRow=0
    \foreach \Row in {#2} {
       \tableauCol=1
       \foreach\k in \Row {
          \draw(\the\tableauCol,\the\tableauRow)+(-.5,-.5)rectangle++(.5,.5);
          \draw(\the\tableauCol,\the\tableauRow)node{\k};
          \global\advance\tableauCol by 1
       }
       \global\advance\tableauRow by -1
    }
  \end{tikzpicture}
}

\DeclarePairedDelimiterX{\set}[1]{\{}{\}}{\setargs{#1}}
\NewDocumentCommand{\setargs}{>{\SplitArgument{1}{|}}m}
{\setargsaux#1}
\NewDocumentCommand{\setargsaux}{mm}
  {\IfNoValueTF{#2}{#1} {#1\,\delimsize|\,\mathopen{}#2}}%{#1\:;\:#2}

\def\({\big(}
\def\){\big)}
\let\<=\langle
\let\>=\rangle
\def\Prod{\displaystyle\prod}

\let\gedom\trianglerighteq
\def\bijection{\overset{\sim}{\longrightarrow}}

\def\Xrightarrow#1{\xRightarrow{\hspace*{1mm}#1\hspace*{1.6mm}}}
\def\xxrightarrow#1{\xrightarrow{\hspace*{1mm}#1\hspace*{1.6mm}}}

\def\half{{\tfrac12}}
\def\q{q^{\frac12}}

\def\N{\mathbb N}
\def\Z{\mathbb Z}

\def\A{{\mathscr A}}
\def\MAn{\mathscr{A}_q(n)}

\def\CC{{\mathscr C}}

\newcommand\Cmin[1][C]{{#1}_{\text{min}}}
\def\H{{\mathscr H}}
\newcommand\Hpm[1][\pm]{\H_q^{#1}}
\newcommand\Hn[1][q]{\H_{#1}(\Sym_n)}
\newcommand\An[1][q]{\H_{#1}(\Alt_n)}
\def\L{{\mathscr L}}
\DeclareMathAlphabet{\mathpzc}{OT1}{pzc}{m}{it}
\def\SS{{\mathpzc S}}
\def\T{{\mathpzc T}}

\newcommand\Em[1][m]{\mathcal E_{#1}}
\def\Alt{\mathfrak A}
\def\Sym{\mathfrak S}
\newcommand\Parts[1][n]{\mathcal{P}_{#1}}
\newcommand\AParts[1][n]{\Parts[#1]^{\text{alt}}}

\def\F{\mathbb F}

\def\Zcal{\mathcal{Z}}
\def\bar{\rule[1.4ex]{.6em}{.1ex}}
\def\s{\mathfrak s}
\def\th{{\text{th}}}
\def\t{\mathfrak t}

\def\v{\mathfrak v}

\let\len=\ell

\def\map#1#2{\,{:}\,#1\!\longrightarrow\!#2}
\def\isomap#1#2{\,{:}\,#1\!\overset{\cong}{\longrightarrow}\!#2}

\def\diag(#1){\llbracket#1\rrbracket}
\def\Diag{\mathop{\rm Diag}}

\newcommand\Mod{{-}\mathop{\rm Mod}}

\DeclareMathOperator{\Span}{span}
\DeclareMathOperator{\Ind}{Ind}
\DeclareMathOperator{\Res}{Res}
\newcommand\TStd[1][\lambda]{\overrightarrow{\text{\rm Std}}(#1)}
\newcommand\Std{\mathop{\rm Std}}

\def\Email#1{\email{\href{mailto:#1}{#1}}}

\begin{document}\bibliographystyle{andrew}
\title[The irreducible characters of the alternating Hecke algebras]
      {The irreducible characters of the\\alternating Hecke algebras}
\author{Andrew Mathas}
\address{School of Mathematics and Statistics, %
         University of Sydney, NSW 2006, Australia.}
\Email{andrew.mathas@sydney.edu.au}

\author{Leah Neves}
\address{School of Mathematics and Statistics, %
         University of Sydney, NSW 2006, Australia.}
\Email{L.Neves@maths.usyd.edu.au}

\begin{abstract}
This paper computes the irreducible characters of the alternating Hecke
algebras, which are deformations of the group algebras of the alternating
groups. More precisely, we compute the values of the irreducible characters of
the semisimple alternating Hecke algebras on a set of elements indexed by
minimal length conjugacy class representatives and we show that these character
values determine the irreducible characters completely. As an application we
determine a splitting field for the alternating Hecke algebras in the semisimple
case.
\end{abstract}
\maketitle

\section*{Introduction}
Mitsuhashi~\cite{Mitsuhashi:A} introduced a deformation, or $q$--analogue,
of the group algebras of the alternating groups. He defined these algebras
by generators and relations. These algebras can also be realized as the
fixed point subalgebra of the Iwahori--Hecke algebra of the corresponding
symmetric group. In the semisimple case, Mitsuhashi showed that all of the
irreducible representations of these algebras can be obtained by Clifford theory
from the irreducible representations of the Hecke algebras of type~$A$.

The aim of this paper is to explicitly compute the irreducible characters of the
alternating Hecke algebras in the semisimple case. As with the alternating
groups, most of the irreducible representations of the Iwahori--Hecke algebra of
type $A$ restrict to give irreducible representations of the alternating Hecke
algebra. The representations that do not remain irreducible on restriction
split as a direct sum of two non-isomorphic irreducible representations so it
only these characters that we need to consider. We compute these characters by
extending ideas of Headley~\cite{Headley}.

One striking feature of the Iwahori--Hecke algebras is that their
characters are determined by their values on any set of standard basis
elements indexed by a set of minimal length conjugacy class
representatives. This was first proved in type~$A$ by
Starkey~\cite{Starkey} (and later rediscovered by Ram~\cite{Ram}), and
then proved for all finite Coxeter groups by Geck and
Pfeiffer~\cite{GeckPfeiffer:irred}. In the final section of this paper
we show that an analogue the Geck--Pfeiffer theorem holds for the
alternating Hecke algebras. Unlike in the Geck-Pfeiffer theorem, the
characters of the alternating Hecke algebras are not constant on minimal
length conjugacy class representations and, by necessity, our proof
relies on a brute force calculation with characters.

The outline of this paper is as follows. In section~1 we define the
alternating Hecke algebra and prove some basic results about it.
Section~2 constructs the irreducible representations of the alternating
Hecke algebras and reduces the calculation of their characters to what
is essentially a calculation in the Iwahori--Hecke algebra of the
symmetric group. Section~3 computes the irreducible characters for those
representations of the Iwahori--Hecke algebras that split on
restriction.  In section~4 we prove a weak analogue of the
Geck--Pfeiffer theorem for the alternating Hecke algebras and, by way of
examples, show that a stronger result is not possible for the bases we
consider. As an application, we determine a splitting field for the
semisimple alternating Hecke algebras.

\section{The alternating Hecke algebra}\label{S:Hecke}
In this section we define the alternating Hecke and construct a basis for it.

Through out this section let $R$ be a commutative ring with one such that
$2$ is a unit in~$R$. Let~$q$ be an invertible element of~$R$ such that
$q\ne1$.

Fix a positive integer $n\ge1$ and let $\Sym_n$ be the symmetric group of
degree~$n$.

\begin{Definition}\label{D:Hecke}
The \textbf{Iwahori--Hecke $R$-algebra of type $A$}, with parameter
$q$, is the unital associative $R$-algebra $\Hn=\H_{R,q}(\Sym_n)$ with generators
$T_1,\dots,T_{n-1}$ and relations
\begin{xalignat*}{2}\label{H relations}
  (T_i-q)(T_i+q^{-1})&=0,                &&\text{for } 1\le i<n,\\
  T_jT_i&=T_jT_i,                   &&\text{for } 1\le i<j-1<n-2,\\
  T_iT_{i+1}T_i&=T_{i+1}T_iT_{i+1}, &&\text{for } 1\le i<n-2.
\end{xalignat*}
\end{Definition}

Define $s_i=(i,i+1)\in\Sym_n$ and $T_{s_i}=T_i$, for $i=1,\dots,n-1$. Then
$S=\{s_1,\dots,s_{n-1}\}$ is the standard set of Coxeter generators for
$\Sym_n$. Every element $w\in\Sym_n$ can be written in the form
$w=s_{i_1}\dots s_{i_k}$, for some~$i_j$ with $1\le i_j\le k$ for
$j=1,\dots,k$. Such an expression is \textbf{reduced} if $k$ is minimal, in
which case we write $\ell(w)=k$ and say that $w$ has \textbf{length}~$k$.

Suppose that $x=s_{i_1}\dots s_{i_k}$ is a reduced expression for $n$.  Set
$T_x=T_{i_1}\dots T_{i_k}$ and $\len(x)=k$. Then, as is well--known,
$\len\map{\Sym_n}\N$ is the usual Coxeter length function on $x$ and $T_x$ is
independent of the choice of $i_1,\dots,i_k$.  Moreover, $\set{T_x|x\in\Sym_n}$
is a basis of $\Hn$. Proofs of all of these facts can be found, for example, in
\cite[Chapter~1]{M:Ulect}.

\begin{Remark}\label{R:Scaling}
  We have rescaled the generators of $\Hn$ when compared with Iwahori's
  original definition~\cite{Iwahori:Hecke} of $\Hn$. As a result, when comparing
  our results with the literature if is necessary to replace $q$ with
  $\q$ and $T_w$ with $q^{-\half\ell(w)}T_w$, for $w\in\Sym_n$. The
  algebra
  $\Hn$  is cellular by \cite[Theorem~3.20]{M:Ulect}, so any field is a
  splitting field for $\Hn$.  Therefore, adjoining a square root of~$q$ to the
  ground ring does not affect the representation theory of~$\Hn$, however,
  adjoining a square root does affect the representation theory of the
  alternating Hecke that we introduce below.  Later we place further
  conditions on the ground ring to ensure that we are working over a
  splitting field for the semisimple alternating Hecke algebras.
\end{Remark}

By the first relation in \cref{D:Hecke}, $T_s$ is invertible with
$T_s^{-1}=T_s-q+q^{-1}$, for $s\in S$.  Consequently, $T_x$ is invertible, for
$x\in\Sym_n$.  For $x\in\Sym_n$ set $\varepsilon_x=(-1)^{\ell(x)}$.

\begin{Definition}\label{D:AltHecke}\hspace*{20mm}
  \begin{enumerate}
  \item Let $\#$ be the unique $R$-linear automorphism of $\Hn$ such that
    $T_x^\#=\varepsilon_xT_{x^{-1}}^{-1}$, for all $x\in\Sym_n$.
  \item The \textbf{alternating Hecke algebra}, with parameter $q$, is the
  $\Hn$-fixed point subalgebra
  \[\An=\set{a\in\Hn|a^\#=a} \]
  \end{enumerate}
\end{Definition}

Iwahori~\cite[Theorem~5.4]{Iwahori:Hecke} attributes the involution $\#$ to
Goldman.

The algebra $\An$ was first considered by Mitsuhashi\cite{Mitsuhashi:A}
who defined it by generators and relations. We show in
\cref{P:Generators} below that \cref{D:AltHecke} agrees with
Mitsuhashi's definition.

Recall that the \textbf{alternating group} $\Alt_n$ is subgroup of $\Sym_n$
consisting of even permutations. Explicitly,
$\Alt_n=\set{z\in\Sym_n|\len(z)\equiv0\pmod 2}$. Note that if $q^2=1$ then
$\Hn[\pm1]\cong R\Sym_n$ and $\An[\pm1]\cong R\Alt_n$ since, in
this case, $T_z^\#=\varepsilon_xz$ for all $z\in\Sym_n$.

% If $h\in\Hn$ and $z\in\Sym_n$ let $(h:z)$ be the coefficient of $T_z$
% when $h$ is written with respect to the $\{T_z\}$ basis of $\Hn$. Thus,
% $h=\sum_z(h:z)T_z$.

Let $\le$ be the \textbf{Bruhat order} on~$\Sym_n$. Thus, if $z=s_{i_1}\dots
s_{i_k}$ is a reduced expression for~$z$ then $y\le z$ if and only if
$y=(j_1,j_1+1)\dots(j_l,j_l+1)$, where $j_a=i_{f(a)}$ for some
monotonically increasing function $f\map{\{1,\dots,l\}}\{1,\dots,k\}$. If
$y\le z$ and $y\ne z$ we write $y<z$ or $z>y$.
Write $x\prec y$ if $x<y$ and $\ell(y)\not\equiv\ell(y)\pmod 2$ and
$x\preceq y$ if $x\prec y$ or~$x=y$.

We start by describing a new basis for $\Hn$ that is compatible with the
alternating Hecke algebra $\An$. For $z\in\Sym_n$ set
\begin{equation}\label{E:Az}
A_z=\tfrac12(T_z+\varepsilon_xT_z^\#).
\end{equation}
For example, if $s\in S$ then $A_s=T_s-\half(q-q^{-1})$.

Set $\Hpm=\set{h\in\Hn|h^\#=\pm h}$. By definition $\Hpm[+]$
and $\Hpm[-]$ are $R$-submodules of $\Hn$, however, note that
$\An=\Hn^+$ is a subalgebra. Abusing notation,
write $\Hpm[\varepsilon_z]=\Hpm$ for the corresponding choice of sign.

\begin{Lemma}\label{L:AxProperties}
    Suppose that $z\in\Sym_n$. Then $A_z=T_z+\sum_{y<z}a_{yz}T_y$, for some
    $a_{yz}\in R$. Moreover, $A_z\in\Hpm[\varepsilon_z]$.
\end{Lemma}

\begin{proof} By definition, $A_z^\#=\varepsilon_xA_z$, so
  $A_z\in\Hpm[\varepsilon_z]$. To complete the proof observe that if
  $z=s_{i_1}\dots s_{i_k}$ is a reduced expression for~$z$ then, as is
  well-known,
  \[T_z^\# =(-T_{s_1}+q-q^{-1})\dots(-T_{s_k}+q-q^{-1})
           =\varepsilon_z T_z+\sum_{y<z}r_{yz} T_y,\]
  for some $r_{yz}\in R$. Hence, the second statement follows by setting
  $a_{yz}=\half\varepsilon_zr_{yz}$.
\end{proof}

\begin{Proposition}\label{P:AnBAsis}
  Suppose that $n\ge1$. Then the following hold:
  \begin{enumerate}
    \item $\set{A_z|z\in\Alt_n}$ is a basis of $\An=\Hpm[+]$.
    \item $\set{A_z|z\in\Sym_n\setminus\Alt_n}$ is a basis of $\Hpm[-]$.
    \item As an $R$-module, $\Hn=\Hpm[+]\oplus\Hpm[-]$.
\end{enumerate}
\end{Proposition}

\begin{proof}By \cref{L:AxProperties}, the transition matrix between the
  $T$-basis $\set{T_z|z\in\Sym_n}$ and the basis $\set{A_z|z\in\Sym_n}$ is
  unitriangular, so $\set{A_z|z\in\Sym_z}$ is a basis of~$\Hn$. Applying
  \cref{L:AxProperties} again, $A_z\in\Hpm[\varepsilon_z]$. Hence, all parts
  of the proposition now follow since $\Hpm[+]\cap\Hpm[-]=0$.
\end{proof}

In particular, $\set{A_w|w\in\Alt_n}$ is a basis of $\An$. When $q^2=1$
the basis element $A_w$, for $w\in\Alt_n$ coincides with the basis
element $w$ in the group ring $R\Alt_n$, however, $A_w$ typically
involves many terms $T_y$ with $y\le w$ and $y\notin\Alt_n$. The algebra
$\An$ has another basis $\set{B_z|z\in\Alt_n}$ that also coincides with
the group-basis of $R\Alt_n$ when $q^2=1$ and which does not involve any
terms $T_y$ with $y<w$ and $y\in\Alt_n$

\begin{Proposition}\label{P:BBasis}
  \begin{enumerate}
    \item For each $z\in\Sym_n$ there exists a unique element
    $B_z\in\Hpm[\varepsilon_z]$ such that
    \[ B_z = T_z + \sum_{y\prec z}b_{yz}T_y, \]
    for some $b_{yz}\in R$.
    \item $\set{B_z|z\in\Alt_n}$ is a basis of $\An=\Hpm[+]$.
    \item $\set{B_z|z\in\Sym_n\setminus\Alt_n}$ is a basis of $\Hpm[-]$.
\end{enumerate}
\end{Proposition}

\begin{proof} Parts~(b) and~(c) follow directly from part~(a) and
  \cref{P:AnBAsis}, so it remains to prove part~(a). First suppose that there
  exist two elements $B_z$ and $B_z'$in $\Hpm[\varepsilon_z]$ that are of the
  required form.  By construction, $B_z-B_z'$ is a $ R$-linear combination of
  terms $T_y$ with $y<z$ and $\varepsilon_y=-\varepsilon_z$. On the other hand,
  $B_z-B_z'\in\Hpm[\varepsilon_z]$ so that $B_z-B_z'$ is a $ R$-linear
  combination of terms $A_y$, with $\varepsilon_y=-\varepsilon_z$. The only way
  that both of these constraints are possible is if $B_z-B_z'=0$. Therefore,
  $B_z$ is uniquely determined as claimed.

  To prove existence we argue by induction on $\len(z)$. If $\len(z)\le1$ then
  $B_z=A_z$ so there is nothing to prove.  If $\len(z)>1$ then the
  element
  \[ B_z=A_z-\sum_{\substack{y<z\\\ell(y)-\ell(z)\in 2\Z}} a_{yz}B_y \]
  has the required properties, where the coefficients $a_{yz}\in R$ are
  given by \cref{L:AxProperties}.
\end{proof}

Set $b_{zz}=1$. The next result describes how $B_r$ acts on the $B$-basis.

\begin{Corollary}\label{C:BAction}
  Suppose that $z\in \Sym_n$ and $1\le r<n$. Then
  \[ B_r B_z=\delta_{rz>z}B_{rz}
      + \half(q-q^{-1})\Bigl(\sum_{\substack{y\prec z\\ry<y}} b_{yz}B_y
      - \sum_{\substack{y\prec z\\ry>y}} b_{yz}B_y\Bigr),
  \]
  where $\delta_{rz>z}=1$ if $rz>z$ and $\delta_{rz>z}=0$ if $rz<z$.
\end{Corollary}

\begin{proof}
  Since $B_r=A_r=T_r-\half(q-q^{-1})$, we compute
  \begin{align*}
    B_rB_z &=\bigl(T_r-\half(q-q^{-1})\bigr)\sum_{y\preceq z}b_{yz}T_y\\
    &=\sum_{\substack{y\preceq z\\ry<y}}
           b_{yz}\bigl((q-q^{-1})T_y+T_{ry}\bigr)
           +\sum_{\substack{y\preceq z\\ry>y}}b_{yz}T_{ry}
           -\half(q-q^{-1})\sum_{y\preceq z}b_{yz}T_y.
  \end{align*}
  If $B_z\in\Hpm[\pm]$ then $B_rB_z\in\Hpm[\mp]$ so that
  $B_rB_z=\sum_y c_yB_y$, where $c_y\ne0$ only if
  $\ell(y)\not\equiv\ell(z)\pmod 2$. By \cref{P:BBasis}, if
  $\ell(y)\not\equiv\ell(z)$ then $c_y$ is equal
  to the coefficient of~$T_y$ in the last displayed equation, so the
  result follows.
  %That is,
  %\[
  %    c_y = \begin{cases}
  %            b_{yz}(q-q^{-1})(1-\half),&\text{if }ry<y\prec z\\
  %            -\half b_{yz}(q-q^{-1}),&\text{if }ry>y\prec z,\\
  %            1,&\textit{if }y=rz>z.
  %          \end{cases}
  %\]
  %giving the result.
\end{proof}

\begin{Remark}
  The argument of \cref{C:BAction} shows that if $rz<z$ then
  $b_{rz,z}=\half(q-q^{-1})$.
\end{Remark}

The basis $\set{B_z|z\in\Alt_n}$ of $\An$ was independently discovered by
Shoji~\cite{Shoji:invol} using a construction that is reminiscent of
that for the Kazhdan-Lusztig basis of~$\Hn$. Shoji's argument gives more
information about the coefficients $b_{yz}$, for $z,y\in\Sym_n$.  Shoji
also showed that each $B_z$ is invariant under Kazhdan and Lusztig's bar
involution~\cite{KL}. This fact will be useful for us below, so we now
give a proof of this.

Assume for the next results that $R=\Zcal$, where $\Zcal=\Z[\half,
q,q^{-1}]$ and~$q$ is an indeterminate over~$\Z$.
The \text{bar involution}~$\bar\map\Zcal\Zcal$ is the unique ring
involution of~$\Zcal$ such that $\overline{q}=q^{-1}$. A
\textbf{semilinear} involution on~$\Hn$ is a $\Z$-linear ring involution
$\iota\map\Hn\Hn$ such that $\iota(q)=q^{-1}$. Inspecting
\cref{D:Hecke}, there are two semilinear involutions of the  Hecke
algebra $\Hn$ given by
\[ \beta(T_z)=T_{z^{-1}}^{-1}\qquad\text{and}\qquad\epsilon(T_z)=\varepsilon_zT_z, \]
for all $z\in\Sym_n$. (The involution $\beta$ is the Kazhdan-Lusztig bar
involution.)

\begin{Corollary}[Shoji~\cite{Shoji:invol}]\label{C:KLInvariance}
  Suppose that $z\in\Sym_n$. Then $\epsilon(B_z)=\varepsilon_zB_z$ and
  $\beta(B_z)=B_z$.
\end{Corollary}

\begin{proof}By \cref{P:BBasis},
  $\epsilon(B_z)=\varepsilon_z T_z+\sum_{y<z}\varepsilon_y\overline{b_{yz}}T_y$,
  where $b_{yz}\ne0$ only if $\varepsilon_y=-\varepsilon_z$. Therefore,
  $\epsilon(B_z)=\varepsilon_z B_z$ by the uniqueness clause of
  \cref{P:BBasis} (consequently,
  $\overline{b_{yz}}=\varepsilon_z\varepsilon_y b_{yz}=-b_{yz}$). By comparing
  the effect of these involutions on the generators of~$\Hn$ shows
  that these two involutions commute and that $\beta=\epsilon\circ\#$.
  Therefore, $\beta(B_z)=\epsilon(B_z^\#)=\varepsilon_x\epsilon(B_x)=B_z$ as
  required.
\end{proof}

\section{Mitsuhashi's generators and relations}\label{S:Mitsuhashi}
  Mitsuhashi~\cite{Mitsuhashi:A} defined an alternating Hecke algebra
  using generators and relations. In this section we show that
  Mitsuhashi's algebra is isomorphic to~$\An$. We assume that $R$ is a
  commutative ring with~$1$ such that $2$ and $q+q^{-1}$ are both
  invertible in~$R$.

  Following Mitsuhashi~\cite[Definition~4.1]{Mitsuhashi:A}, but taking
  \cref{R:Scaling} into account, define the \textbf{Mitsuhashi alternating Hecke
  algebra} $\MAn$ to be the unital associative $R$-algebra
  with generators $A_1,\dots,A_{n-1}$ and relations
  \begin{align*}
    A_1=1, \qquad (A_iA_j)^2=1\quad\text{and}\quad
    \Bigl((A_{i-1}A_i)^2+\Bigl(\frac{q-q^{-1}}{q+q^{-1}}\Bigr)^2
    (A_{i-1}A_i-1)\Bigr)A_{i-1}A_i=1,
  \end{align*}
  where $2\le i,j<n$ and $|i-j|\ne1$. (This is a streamlined version of
  Mitsuhashi's presentation.)

  Mitsuhashi obtained this presentation by first giving a new presentation
  for $\Hn$. For $1\le i<n$ define
  \[E_i=\frac{(2T_i-q+q^{-1})}{(q+q^{-1})}.\]
  We record some elementary properties of these elements (except for the
  first claim, these can be found in~\cite{Mitsuhashi:A}).

  \begin{Lemma}\label{L:EiPropertiess}
      Suppose that $1\le i<n$. Then $E_i^\#=-E_i$, $E_i^2=1$ and
      $T_i=\frac{(q+q^{-1})}2E_i+\half(q-q^{-1})$.
  \end{Lemma}

  In particular, since $T_i=\half((q+q^{-1})E_i+q-q^{-1})$, it follows
  that $\Hn$ is generated as an $R$-algebra by the elements
  $\set{E_i|1\le i<n}$.

\begin{Proposition}\label{P:Generators}
  Suppose that $2$ and $q+q^{-1}$ are invertible in $R$. Then
  $\An\cong\MAn$ as $R$-algebras.
\end{Proposition}

\begin{proof}
  By \cref{L:EiPropertiess}, the set
  $\set{E_1E_i|2\le i<n}$ generates $\An$ as an $R$-algebra since
  $\set{E_i|1\le i<n}$ generates~$\Hn$.  As noted in
  \cite[Lemma~4.2]{Mitsuhashi:A} is follows easily from the relations
  in~$\Hn$ that $E_iE_j=E_jE_i$ if $|i-j|\ne1$ and that
  \[ \Bigl((E_{i-1}E_i)^2+\Bigl(\frac{q-q^{-1}}{q+q^{-1}}\Bigr)^2
                        (E_{i-1}E_i-1)\Bigr)E_{i-1}E_i=1
  \]
  for $2<i<n$. Hence, there is a unique algebra homomorphism
  $\theta\map\MAn\An$ such that $\theta(A_i)=E_1E_i$, for $1\le i<n$. By
  \cite[Theorem~4.7]{Mitsuhashi:A}, $\MAn$ is free as an $R$-module of
  rank $\half n!$, so $\theta$ is an isomorphism by \cref{P:AnBAsis}.
\end{proof}

As a consequence of the proof of \cref{P:Generators}, if $2$ and
$q+q^{-1}$ are both invertible then the $R$-module decomposition of
$\Hn$ from \cref{P:AnBAsis}(c) becomes an $\An$-module decomposition.

\begin{Corollary}\label{super}
  Suppose that $2$ and $q+q^{-1}$ are invertible in $R$. Then, as an
  $\An$-module
  \[ \Hn=\An\oplus E_1\An. \]
  In particular, $\Hn$ is free as an $\An$--module.
\end{Corollary}

Let $M$ be an $\Hn$--module. Define $M^\#=\set{m'|m\in M}$ to be the
$\Hn$--module that is isomorphic to $M$ as a vector space but where the
$\Hn$--action twisted by $\#$. More explicitly, if $m\in M$ and
$h\in\Hn$ then $m'h=(mh^\#)'$. Observe that $M\cong M^\#$ as
$\An$--modules since $a^\#=a$, for all $a\in\An$.

By \cref{super} and a standard application of Clifford
theory (see, for example, \cite[appendix]{RamRamagge}), we have:

\begin{Proposition}\label{P:CliffordTheory}
  Suppose that $\F$ is an algebraically closed field in which both
  $2$ and $q+q^{-1}$ are invertible. Let~$M$ be an irreducible
  $\Hn$--module. Then $M$ is irreducible as an $\An$--module if and only
  if $M\not\cong M^\#$ as $\Hn$--modules. Moreover, $M\cong M^\#$ as
  $\Hn$--modules if and only if $M\cong M^+\oplus M^-$ as
  $\An$--modules, where $M^+$ and $M^-$ are non--isomorphic irreducible
  $\An$-modules.
\end{Proposition}

\begin{proof}
  In the notation of \cite[Appendix]{RamRamagge}, take $A=\An$. Then
  $\Hn\cong\An\ltimes\Z/2\Z$. The condition that $M\cong M^\#$
  is equivalent to saying that the inertia group of $M$ is equal to~$\Z/2\Z$,
  so the Proposition is a special case of \cite[Theorem~A.6]{RamRamagge}.
\end{proof}

\section{Conjugacy class representatives}\label{S:conjugacy}
Before constructing the irreducible representations of $\An$ we
recall the labeling of the conjugacy classes of $\Sym_n$ and $\Alt_n$.
The full details can be found, for example, in \cite[Chapter~1]{JK}.

A \textbf{composition} of $n$ is a sequence $\mu=(\mu_1,\dots,\mu_k)$ of
positive integers that sum to~$n$. (For convenience, we set $\mu_0=0$
in many places below.) A composition $\mu$ is a \textbf{partition} if its
parts $\mu_1,\dots,\mu_k$ are in weakly decreasing order.  The
\textbf{length} of a partition or composition is $\ell(\mu)=k$ if $\mu$
has~$k$ parts. Let $\Parts$ be the set of partitions of~$n$.

We identify a partition $\mu$ with its \textbf{diagram}:
\[ \diag(\mu)=\set{(r,c)|1\le c\le\mu_r\text{ for }r\ge1}. \]
In this way, we talk about the rows and columns of (the diagram of) $\lambda$.
If $\lambda$ and $\mu$ are partitions such that $\lambda_i\le\mu_i$, for
$i\ge1$, then $\lambda$ is a \textbf{contained} in $\mu$, $\mu/\lambda$ is a
\textbf{skew partitions}, and we write $\lambda\subseteq\mu$ and
$\lambda\subset\mu$ if $\lambda\ne\mu$.  If $\lambda\subseteq\mu$ then the set
\[\diag(\mu/\lambda)=\diag(\mu)/\diag(\lambda)
        =\set{x\in\diag(\mu)|x\notin\diag(\lambda)}\]
is the \textbf{skew diagram} of shape $\mu/\lambda$.

If $\lambda$ is a partition then its \textbf{conjugate} is the partition
$\lambda'=(\lambda_1',\lambda_2',\dots)$ where $\lambda'_i$ is the number of
entries in column~$i$ of $\diag(\lambda)$. A partition $\lambda$ if
\textbf{self-conjugate}, or \textbf{symmetric}, if~$\lambda=\lambda'$.

As is well-known, and is straightforward to prove, two elements of the symmetric group
are conjugate if and only if they have the same cycle type. Consequently,
the conjugacy classes of~$\Sym_n$ are indexed by the partitions
of~$n$.  If $\kappa=(\kappa_1,\dots,\kappa_k)$ is a composition of~$n$ let
$C_\kappa$ be the elements of $\Sym_n$ of cycle type~$\kappa$.  Then
$\Sym_n=\coprod_{\kappa\in\Parts} C_\kappa$ is the decomposition of $\Sym_n$ into
the disjoint union of its conjugacy classes. If
$\kappa=(\kappa_1,\dots,\kappa_d)$ let
\[w_\kappa=(1,2,\dots,\kappa_1)(\kappa_1+1,\kappa_1+2,\dots,\kappa_1+\kappa_2)
    \dots(\kappa_1+\dots+\kappa_{d-1}+1,\dots,\kappa_1+\dots+\kappa_d).
\]
Then $w_\kappa$ is an element of minimal length in $C_\kappa$. In
particular, $\set{w_\kappa|\kappa\in\Parts}$ is a complete set of
conjugacy class representatives for $\Sym_n$.

Now consider the alternating group $\Alt_n$. If
$\kappa=(\kappa_1,\dots,\kappa_k)$ is a partition of~$n$ then
$w_\kappa\in\Alt_n$ if and only if
$\ell(w_\kappa)=\sum_{i=1}^l(\kappa_i-1)=n-\ell(\kappa)$ is even. Equivalently,
$w_\kappa\in\Alt_n$ if~$\kappa$ has an even number of non--zero even parts. Let
$\AParts=\set{\kappa\in\Parts|n\equiv\ell(\kappa)\pmod 2}$ be this set of
partitions so that $w_\kappa\in\Alt_n$ if and only if $\kappa\in\AParts$.

Fix $\kappa\in\AParts$. If $\kappa$ contains a repeated part, or a non--zero
even part, then~$w_\kappa$ commutes with an element of
odd order, which implies that $C_\kappa$ is a single conjugacy class of
$\Alt_n$.  As a notational sleight of hand, set $w_\kappa^+=w_\kappa=w_\kappa^-$
in this case.  On the other hand, if $n>1$ and all of the parts of
$\kappa\in\AParts$ are all odd and distinct then $C_\kappa$ is the disjoint
union of two conjugacy classes in~$\A_n$. Set $w_\kappa^+=w_\kappa$ and
$w_\kappa^-=s_rw_\kappa s_r$, where $r=\kappa_1+\dots+\kappa_{d-1}+1$
and~$d$ is minimal such that $\kappa_d>1$. Then $w_\kappa^-\in\Alt_n$
and $\ell(w_\kappa^-)=\ell(w_\kappa)$ with $w_\kappa^\pm$ being
conjugate in~$\Sym_n$ but \textit{not} conjugate in~$\Alt_n$. In view of
\cite[Lemma~1.2.10]{JK} and \cite[Example~3.1.16]{GeckPfeiffer:book},
$\set{w_\kappa^\pm|\kappa\in\AParts}$ is a complete set of minimal
length conjugacy class representatives for~$\Alt_n$, for $n>1$.

The description of the conjugacy classes of $\Alt_n$ is not quite what we expect
because, in the semisimple case, it is well-known that the irreducible
representations of~$\Sym_n$ that are indexed by the self-conjugate partitions
split  on restriction to~$\Alt_n$ (see \cref{P:Splitting}). The combinatorial
connection between the representations and the conjugacy classes
that split is given by the following definition.

\begin{Definition}\label{D:hlambda}
  Suppose that $\lambda=\lambda'$.
  Let $h(\lambda)=(h_1,h_2,\dots,h_{d(\lambda)})$, where
  $d(\lambda)=\max\set{i|\lambda_i\ge i}$ is the length of the
  diagonal in $\diag(\lambda)$ and
  $h_i=\lambda_i+\lambda_i'-2i+1$, for $1\le i\le d(\lambda)$.
\end{Definition}

By definition, all of the parts of $h(\lambda)$ are odd and distinct.
Pictorially, $h_i$ is the length of the~$i^\th$ diagonal hook
$H_i(\lambda)=\set{(i,j),(j,i)|i\le j\le\lambda_i}$ in the diagram of $\lambda$.
For example, if $\lambda=(k+1,1^k)$ is a hook partition then $h(\lambda)=(2k+1)$
and $d(\lambda)=1$.

\section{The irreducible representations of $\An$}\label{S:IrreducibleReps}

This section gives an explicit description of the irreducible
representations of $\An$ when $\An$ is semisimple. These representations
have already been described by Mitsuhashi~\cite{Mitsuhashi:A}, however, the
construction that we give of the irreducible $\An$-modules is different to his.
In the next section we use these modules to compute the irreducible
characters of~$\An$.

As we are interested in computing characters we want to work over a field that
is large enough to ensure that $\An$ is a split semisimple
algebra. As we will see, this requires that our ground field contains certain
square roots. To describe these, for any integer $k$ define the
\textbf{$q$--integer} $[k]=[k]_q$ by
\[[k]=\begin{cases}\phantom{-}
  q+q^3+\dots+q^{2k-1},&\text{if }k\ge0,\\
  -q-q^{-3}-\dots-q^{2k+1},&\text{if }k<0.
\end{cases}\]
In particular,  $[0]=0$, $[1]=q$ and $[-1]=-q^{-1}$. If $q^2\ne1$ then
$[k]=(q^{2k}-1)/(q-q^{-1})$, whereas if $q^2=1$ then $[k]=k$ for all
$k\in Z$. By definition, if $k>0$ then $[-k]=-q^{-2k}[k]$, for $k\in\Z$.

\begin{Definition}\label{D:SquareRoots}
  Suppose that $\F$ is a field and that $q\in\F$. We assume that the
  elements $-1$, $[k]$, for $1\le k\le n$, are non-zero and have
  square-roots in~$F$.

  Fix a choice of square roots $\sqrt{-1}$
  and $\sqrt{[k]}$ in~$\F$, for $1\le k\le n$, and set
  $\sqrt{[-k]}=q^{-k}\sqrt{-1}\sqrt{[k]}\in\F$.
\end{Definition}

For example, if $F$ is any field and $q$ is in indeterminate over~$F$
then we could take $\F=F(\sqrt{-1},\sqrt{[1]},\dots,\sqrt{[n]})$, which
is a subfield of~$F((q))$.

For the remainder of this paper we work over the field $\F$ and consider
the $\F$-algebras $\Hn=\H_{\F,q}(\Sym_n)$ and $\An=\H_{\F,q}(\Alt_n)$.
\cref{D:SquareRoots} implies that $[1][2]\dots[n]\ne0$ in~$\F$.  By
results going back to Hoefsmit~\cite{Hoefsmit}, this condition implies
that $\Hn$ is split semisimple. In fact, since $\Hn$ is a cellular
algebra~\cite{M:Ulect}, any field is a splitting field for~$\Hn$. By
\cref{C:SplittingField} below,~$\F$ is a splitting field for $\An$.

If $\lambda$ is a partition (or a skew partition) then a
\textbf{$\lambda$-tableau} is an injective map $\t\map{\diag(\lambda)}\Z$.  We
think of a $\lambda$-tableau as a labeling of the diagram of $\lambda$
and we say that $i$ appears in row~$r$ and column~$c$ of $\t$ if
$\t(r,c)=i$.

A \textbf{standard $\lambda$-tableau} is a bijection
$\t\map{\diag(\lambda)}\{1,\dots,n\}$ such that the entries of this
tableau increase along rows and down columns. More precisely, $\t$ is
standard if $\t(r,c)<\t(s,d)$ whenever $(r,c)$ and $(s,d)$ are distinct
elements of $\diag(\lambda)$ with $r\le s$ and $c\le d$. Let
$\Std(\lambda)$ be the set of standard $\lambda$-tableaux and let
$\Std(\Parts)=\bigcup_{\lambda\in\Parts}\Std(\lambda)$.  If
$\t\in\Std(\lambda)$ then the \textbf{conjugate tableau} $\t'$ is the
standard $\lambda'$-tableau $\t'$ such that $\t'(c,r)=\t(r,c)$, for
$(c,r)\in\diag(\lambda')$.

Let $\t$ be a $\lambda$-tableau and $i$ an integer with $1\le i\le n$.  Then
$\t(r,c)=i$ for some $(r,c)\in\diag(\lambda)$.  The \textbf{content} of $i$ in
$\t$ is the integer $c_\t(i)=c-r$. Observe that if $\t$ is standard then
$c_\t(i)>0$ if $i+1$ appears in a later column of $\t$ than $i$ and, otherwise,
$c_\t(i)<0$.  If $i<n$ then the \textbf{axial distance} from $i$ to $i+1$ is
$\rho_\t(i)=c_\t(i)-c_\t(i+1)$.

To compute the irreducible characters of $\An$ we need an explicit description
of its representations. To this end, for $k\in\Z$ set
\begin{align*}
    \alpha_k&=\begin{cases}
      \frac{\sqrt{-1}\sqrt{[k+1]}\sqrt{[k-1]}}{[k]},&\text{if }k>0,\\
      -\alpha_{-k},&\text{if }k<0.
    \end{cases}\\
\intertext{If $\t$ is a standard $\lambda$-tableau and $1\le i<n$ then set}
 \alpha_\t(i)&=\begin{cases}
                    \alpha_{\rho_\t(i)},&\text{if }\t s_i\in\Std(\lambda),\\
                    0,&\text{otherwise}.
                \end{cases}
\end{align*}
In particular, if $\t$ is standard then $\alpha_{\t'}(i)=-\alpha_\t(i)$,
for $1\le i<n$. This property of the $\alpha$-coefficients is crucial
below because it ensures that the Specht modules for $\Hn$ that are
indexed by self-conjugate partitions split when they are restricted
to~$\Hn$.

\begin{Proposition}\label{P:seminormal}
Suppose that $\lambda$ is a partition and let
$S(\lambda)$ be the $\F$--vector space with basis
$\set{v_\t|\t\in\Std(\lambda)}$. Then $S(\lambda)$ becomes an irreducible
$\Hn$--module with $\Hn$--action for $i=2,\dots,n-1$ given by
\[ v_\t T_i = \frac{-1}{[\rho_\t(i)]}v_\t+\alpha_\t(i)v_{\t s_i}.  \]
\end{Proposition}

Note that if $\rho_\t(i)=\pm1$ then $\t s_i$ is not standard and
$v_\t T_i = \frac{-1}{[\pm1]}v_\t=\mp q^{\mp1}v_\t$ since $\alpha_\t(i)=0$.

\begin{proof}
  Extending the classical arguments for the symmetric group,
  Hoefsmit~\cite{Hoefsmit} has the first to give a seminormal form for the
  irreducible representations of~$\Hn$. This result can be proved by essentially
  repeating Hoefsmit's argument. Alternatively, it is straightforward to check
  that if $\t$ and $\v=\t s_i$ are both standard tableaux then
  \begin{equation}\label{E:SCNCquad}
    \alpha_\t(i)\alpha_\v(i)
      =\frac{[1+\rho_\t(i)][1+\rho_\v(i)]}{q^2[\rho_\t(i)][\rho_\v(i)]}
  \end{equation}
  since $\rho_\v(i)=-\rho_\t(i)$. Moreover, if $\t$ is standard and $1\le i<n-1$
  then
  \[
    \alpha_\t(i)=\alpha_{\t s_{i+1}s_i}(i+1),\quad
    \alpha_{\t s_i}(i+1)=\alpha_{\t s_{i+1}}(i),\quad\text{and}\quad
    \alpha_{\t s_is_{i+1}}(i)=\alpha_\t(i+1),
  \]
  so that
  $\alpha_\t(i)\alpha_{\t s_i}(i+1)\alpha_{\t s_is_{i+1}}(i)
     =\alpha_\t(i+1)\alpha_{\t s_{i+1}}(i)\alpha_{\t s_{i+1}s_i}(i+1)$.
  Similarly, if $|i-j|>1$ then
  $\alpha_\t(i)\alpha_{\t s_i}(j)=\alpha_\t(j)\alpha_{\t s_{j}}(i)$.
  Therefore, the set of scalars
  $\set{\alpha_\t(i)|1\le i<n\text{ and }\t\in\Std(\Parts)}$ is a
  \textit{$*$-seminormal coefficient system} for $\Hn$ in the sense of
  \cite[\S3]{HuMathas:SeminormalQuiver}. (When comparing the
  results of this paper with \cite{HuMathas:SeminormalQuiver}, recall
  \cref{R:Scaling}, and note that \cref{E:SCNCquad} corresponds to the
  quadratic relation $(T_i-q)(T_i+q^{-1})=0$. This accounts for the
  additional factor of~$q^{-2}$ in \cref{E:SCNCquad} when compared with
  the corresponding formula in~\cite{HuMathas:SeminormalQuiver}.) Hence,
  the proposition is a special case of
  \cite[Theorem~3.13]{HuMathas:SeminormalQuiver}.
\end{proof}

% \begin{Remark}\label{R:QuantumIntegers}
%   For $k=1,\dots,n$ let $L_k=\sum_{j=1}^{k-1}T_{(j,k)}$ the be (rescaled)
%   Jucys-Murphy elements of $\Hn$. When $\Hn$ is semisimple, $\<L_1,\dots,L_n\>$
%   is a maximal commutative subalgebra of~$\Hn$. Using \cref{P:seminormal} and
%   arguing by induction on~$k$, a slightly tedious calculation shows that
%   $v_\t L_k=\qm[c_\t(k)]v_\t$, for $1\le k\le n$ and $\t\in\Std(\Parts)$.
% \end{Remark}

We want to describe the Specht modules as $\An$--modules.  Let
$\tau\map{S(\lambda)}{S(\lambda')}$ be the $\F$--linear map given by
$v_\t\mapsto v_{\t'}$, for all $\t\in\Std(\lambda)$. The map $\tau$
depends on the partition $\lambda$, however, its meaning should always
be clear from context.

\begin{Proposition}\label{twist}
    Suppose that $w\in\Sym_n$. Then $v_\t T_w^\#=\tau(v_{\t'} T_w)$,
    for all $\t\in\Std(\lambda)$.
\end{Proposition}

\begin{proof}
Fix $\t\in\Std(\lambda)$. Suppose first that $\ell(w)=1$, so that $w=s_i$
for some $i$. Let $\rho_\t(i)=c_\t(i+1)-c_\t(i)$ be the axial distance from~$i$
to~$i+1$ in~$\t$. If $|\rho_\t(i)|=1$ then $i$ and $i+1$ are in the same row or in
the same column of $\t$, and the Lemma follows directly from the
definitions. Suppose now that $|\rho_\t(i)|>1$. Recall that
$\alpha_{-k}=-\alpha_k$ if $k\ne0$. Therefore,
\begin{align*}
   v_\t T_i^\#&= v_\t (-T_i+q-q^{-1})
     =\big(\frac{1}{[\rho_\t(i)]}+q-q^{-1}\big)v_\t -\alpha_\t(i)v_{\t s_i}\\
    &=\frac{-1}{[-\rho_\t(i)]}v_\t+\alpha_{\t'(i)}v_{\t s_i}
     =\tau\(v_{\t'}T_i\)
\end{align*}
since $\rho_{\t'}(i)=-\rho_\t(i)$ and $(\t s_i)'=\t' s_i$.
The general case now follows easily by induction on~$\ell(w)$.
\end{proof}

By \cref{P:AnBAsis}(c), $\Hn$ is free as an $\An$--module. Therefore,
the natural induction and restriction are exact functors
\[
   \Ind\map{\An\Mod}\Hn\Mod\quad\text{and}\quad\Res\map{\Hn\Mod}\An\Mod.
\]
between the categories of finitely generated $\Hn$-modules and
$\An$-modules, respectively. If $M$ is an $\Hn$-module we abuse
notation and write $M=\Res M$ for the restriction of~$M$ to
an~$\An$-module.

\begin{Corollary}\label{C:ConjugateIso}
Let $\lambda$ be a partition. Then
$\tau\isomap{S(\lambda)}S(\lambda')$ is an isomorphism of
$\An$--modules.
\end{Corollary}

\begin{proof} By definition $\tau$ is an isomorphism of
$\F$--vector spaces, so we only need to show that $\tau$ commutes with the
$\An$--action. Suppose that $a\in\An$ and let $\t\in\Std(\lambda)$.
Then $\tau(v_{\t'} a)=v_{\t'}a^\#=v_{\t'}a=\tau(v_\t)a$ by the Proposition.
\end{proof}

Let $\chi^\lambda$ be the character of the $\Hn$--module $S(\lambda)$
and let $\chi_\A^\lambda$ be the character of the $\An$-module $\Res S(\lambda)$.

\begin{Corollary}\label{C:Hash}
    Suppose that $w\in\Sym_n$. Then
    $\chi^\lambda(T_w^\#)=\chi^{\lambda'}(T_w)$.
\end{Corollary}

\begin{proof} For $\t'\in\Std(\lambda')$ write
  $v_{\t'} T_w=\sum_\s a_{\t'\s'} v_{\s'}$ for some
  scalars $a_{\t'\s'}=a_{\t'\t'}(w)\in\F$ and where, in the sum, $\s'\in\Std(\lambda')$.
  Then $v_\t T_w^\#=\sum_\s a_{\t'\s'}v_{\s}$ by \cref{twist}.
  Therefore, $\chi^\lambda(T_w^\#)=\sum_\t a_{\t'\t'}=\chi^{\lambda'}(T_w)$.
\end{proof}

Recall from \cref{S:Mitsuhashi} that if $M$ is an $\Hn$-module then
$M^\#$ is the $\Hn$ that is equal to~$M$ as a vector space but with the
$\Hn$-action twisted by $\#$.

\begin{Corollary}\label{C:HashIso}
  Suppose that $\lambda$ is a partition of $n$. Then
  $S(\lambda)^\#\cong S(\lambda')$ as $\Hn$--modules.
\end{Corollary}

Combining \cref{C:ConjugateIso} and \cref{C:HashIso},
$S(\lambda)\cong S(\lambda)^\#\cong S(\lambda')$ as $\An$-modules.
Therefore, a standard application of \textit{Clifford theory}
(\cref{P:CliffordTheory}), implies that $S(\lambda)\cong S(\lambda')$ is
either  irreducible $\An$-module if $\lambda\ne\lambda'$. It is
straightforward  to prove this directly by modifying the argument given
in \cref{P:Splitting} below.

The next result shows that if $\lambda\ne\lambda'$ then the irreducible
character $\chi^\lambda_\A$ for~$\An$ is completely determined by the
corresponding irreducible character $\chi^\lambda$ of~$\Hn$. Recall from
\cref{L:AxProperties} that $\set{A_w|w\in\Alt_n}$ is a basis of $\An$.

\begin{Corollary}\label{C:NonSplitCharacters}
  Suppose that $\lambda$ is a partition of $n$ and that
  $w\in\Alt_n$. Then
  \[ \chi^\lambda_\A(A_w)=\half\(\chi^\lambda(T_w)+\chi^{\lambda'}(T_w)\). \]
\end{Corollary}

It remains to consider the irreducible characters $\chi^{\lambda\pm}_\A$ of~$\An$, where
$\lambda=\lambda'$ is a self--conjugate partition of~$n$. In this
case,~$\tau$ is an $\An$--module automorphism of $S(\lambda)$ by
\cref{C:ConjugateIso}. A straightforward calculation shows that the
eigenspaces of $\tau$ on $S(\lambda)$ are $\An$--modules. Since
$\tau^2=1$, the only possible eigenvalues for the action of~$\tau$
on~$S(\lambda)$ are~$\pm1$.

\begin{Definition}
  Suppose that $\lambda=\lambda'$. Let $S(\lambda)^+$ and $S(\lambda)^-$ be the
  vector subspaces
   \[ S(\lambda)^\pm=\Span\set{v_\t\pm v_{\t'}|\t\in\TStd}, \]
   where $\TStd$ is the set of standard $\lambda$-tableau that have
   $2$ in their first row.
\end{Definition}

\begin{Proposition}\label{P:Splitting}
  Suppose that $\lambda=\lambda'$. Then $S(\lambda)^+$ and $S(\lambda)^-$ are
  irreducible $\An$-submodules of $S(\lambda)$. Moreover,
  $S(\lambda)=S(\lambda)^+\oplus S(\lambda)^-$ as $\An$-modules.
\end{Proposition}

\begin{proof}By inspection, $S(\lambda)^\pm$ is the $\pm1$-eigenspace for~$\tau$
  acting on~$S(\lambda)$. Therefore, $S(\lambda)^\pm$ is an $\An$-submodule of
  $S(\lambda)$ and $S(\lambda)^+\cap S(\lambda)^-=0$. Now, if
  $\t\in\Std(\lambda)$ then either~$2$ is in the first row of~$\t$ and
  $\t\in\TStd$ or $2$ is in the first column of~$\t$ and
  $\t'\in\TStd$. Hence,
  $S(\lambda)=S(\lambda)^+\oplus S(\lambda)^-$. Finally, the two modules
  $S(\lambda)^\pm$ are irreducible $\An$-modules by Clifford theory.
  Rather than using Clifford theory we can
  prove this directly using the fact that if $F_\t\in\Hn$ is the primitive
  idempotent attached to the standard $\lambda$-tableau~$\t$, as in
  \cite[Theorem~3.13]{HuMathas:SeminormalQuiver}, then
  $F_\t+F_\t^\#=F_\t+F_{\t'}\in\An$ since $F_\t^\#=F_{\t'}$ by
  \cite[Lemma~2.7]{BM:AlternatingKLR}.  Consequently, if
  $x=\sum_{\s\in\TStd(\lambda)}r_\s(v_\s\pm v_{\s'})$, for some
  $r_\s\in\F$, is any non-zero element of $S(\lambda)^\pm$ with $r_\t\ne0$
  then
  \[ v_\t\pm v_{\t'}=\frac1{r_\t}x(F_\t+F_{\t'})\in x\An.\]
  Using \cref{P:seminormal} it follows that
  $x\An=S(\lambda)^\pm$, showing that $S(\lambda)^\pm$ is an irreducible
  $\An$-module.
\end{proof}

Let $\succ$ be the lexicographic order on $\Parts$. The results
above show that
\[  \set{S(\lambda)|\lambda\succ\lambda'\text{ for }\lambda\in\Parts}
       \cup\set{S(\lambda)^\pm|\lambda=\lambda'\text{ for }\lambda\in\Parts}
\]
is a complete set of pairwise non-isomorphic absolutely irreducible
$\An$-modules. The classification of the irreducible $\An$-modules in
the semisimple case is due to Mitsuhashi~\cref{P:CliffordTheory}. As we
have an explicit realisation of all of the irreducible representations
of~$\An$ over~$\F$, and because these modules cannot split anymore
according to Clifford theory, we obtain the following.

\begin{Corollary}\label{C:SplittingField}
  The field $\F=F(\sqrt{-1},\sqrt{[1]}, \sqrt{[2]},\dots,\sqrt{[n]})$
  is a splitting field for $\An[\F,q]$.
\end{Corollary}

Let $\chi^{\lambda\pm}_\A$ be the character of the $\An$-module
$S(\lambda)^\pm$, respectively. We want to compute
$\chi^{\lambda\pm}_\A(a)$, for $a\in\An$. The next result shows that we
can reduce the calculation of these characters to what is essentially
a computation inside $\Hn$. Before we can state this result we need
some more notation.

If $h\in\Hn$ then right multiplication by $h$ gives an endomorphism
$\rho_h$ of $S(\lambda)$. Let $h\tau=\rho_h\circ\tau$ be the
endomorphism obtained by composing this map with $\tau$. Then
$\chi^\lambda(h\tau)$ is the trace of $h\tau$ acting on
$S(\lambda)$.

\begin{Proposition}\label{P:SplitCharacters}
  Suppose that $\lambda=\lambda'$ and $a\in\An$.  Then
  \[ \chi^{\lambda\pm}_\A(a)=\half\(\chi^\lambda(a)\pm\chi^\lambda(a\tau)\). \]
\end{Proposition}

\begin{proof}
The direct sum decomposition
$S(\lambda)=S(\lambda)^+\oplus S(\lambda)^-$ implies that
$\chi^\lambda(a)=\chi_\A^{\lambda+}(a)+\chi_\A^{\lambda-}(a)$, for
$a\in\An$. Further,
$\set{v_\t\pm v_{\t'}|\t\in\TStd}$ is a basis of $S(\lambda)^\pm$ by
\cref{P:Splitting}. Computing traces with respect to this basis shows
that $\chi^\lambda(a\tau)=\chi_\A^{\lambda+}(a)-\chi_\A^{\lambda-}(a)$.
Combining these two formulas,
    $\chi^{\lambda\pm}_\A(a)
       =\half\(\chi^\lambda(a)\pm\chi^\lambda(a\tau)\)$
as claimed.
\end{proof}

\begin{Corollary}\label{C:SplitCharacters}
  Suppose that $\lambda=\lambda'$ and that $w\in A_n$. Then
    \[\chi^{\lambda\pm}_\A(A_w)
        =\half(\chi^\lambda(T_w)\pm\chi^\lambda(T_w\tau)).\]
\end{Corollary}

\begin{proof}By \cref{P:SplitCharacters},
\begin{align*}
\chi^{\lambda\pm}_\A(A_w)
   &=\frac14\(\chi^\lambda(T_w+T_w^\#)\pm\chi^\lambda((T_w+T_w^\#)\tau)\)\\
   &=\half\chi^\lambda(T_w)\pm\frac14\(\chi^\lambda(T_w\tau)
               +\chi^\lambda(T_w^{\#}\tau)\),
\end{align*}
by \cref{C:Hash} since
$\chi^\lambda(T_w^\#)=\chi^{\lambda'}(T_w)=\chi^\lambda(T_w)$.
Arguing as in \cref{C:Hash} shows that
$\chi^\lambda(T_w\tau)=\chi^\lambda(T_w^{\#}\tau)$, which
completes the proof.
\end{proof}

As the character values $\chi^\lambda(T_w)$ are known for $w\in\Sym_n$
(see, for example, \cite{GeckPfeiffer:book}), we are reduced to
computing $\chi^\lambda(T_w\tau)$, for $w\in A_n$ when
$\lambda=\lambda'$ is a self--conjugate partition. This is the subject
of next section.

\section{The character values for self--conjugate partitions}
The results of the last section construct the irreducible characters of~$\An$
and shows that these characters are determined by the irreducible characters
of~$\Hn$ except, possibly, for the characters $\chi^\lambda_\A$ when
$\lambda=\lambda'$. When $\lambda=\lambda'$ is a self-conjugate partition
\cref{C:SplitCharacters} reduced the calculation of $\chi^\lambda(A_w)$
to understanding $\chi^\lambda(T_w\tau)$, for $w\in\Alt_n$. This
section computes the character values
$\chi^\lambda_\H(T_{w_\kappa}\tau)$ whenever $\lambda$ is a
self-conjugate partition and $\kappa=(\kappa_1,\dots,\kappa_d)$ is a
composition of~$n$.

The next theorem is one of the main results of this paper. It determines
the irreducible characters~$\chi^{\lambda\pm}_\A$ of~$\An$ when
$\lambda$ is a self-conjugate partition. Recall that we fixed a
field~$\F$ in \cref{D:SquareRoots} and if~$\lambda=\lambda'$ is
self-conjugate then the partition
$h(\lambda)=(h_1,h_2,\dots,h_{d(\lambda)})$ is given in
\cref{D:hlambda}.

\begin{Theorem}\label{T:TauCharacters}
  Suppose that $\lambda=\lambda'$ is a self--conjugate partition of~$n$ and
  let~$h(\lambda)=(h_1,h_2,\dots,h_{d(\lambda)})$. Let $\kappa$ be a
  composition of~$n$. Then
  \[\chi^\lambda(T_{w_\kappa}\tau)
               =\begin{cases}\displaystyle
                 \epsilon_\kappa (-\sqrt{-1})^{\half(n-d(\lambda))}
                 q^{-\half n} \prod_{i=1}^{d(\lambda)}\sqrt{[h_{i}]},&
     \text{if }\vec\kappa=h(\lambda),\\
     0,&\text{otherwise,}
                \end{cases}
  \]
  where $\epsilon_\kappa = (-1)^{\#\set{1\le y<z\le d|\kappa_y<\kappa_z}}$.
\end{Theorem}

In the final section of the paper we show that these character
values determine the irreducible characters of~$\An$ completely.
The reader may check that if $\vec\kappa=h(\lambda)$ then
$\ell(w_\kappa)=(n-d(\lambda))/2$.

Throughout this section we fix a self-conjugate partition
$\lambda=\lambda'$, and a composition $\kappa$, as in
\cref{T:TauCharacters}.

Suppose that $w\in\Sym_n$ and $\t\in\Std(\lambda)$. Define $\gamma_\t(w)\in\F$ to be
the coefficient of $v_{\t'}$ in $v_\t T_w$, so that
\[v_\t T_w=\gamma_\t(w)
v_{\t'}+\sum_{\substack{\s\in\Std(\lambda)\\\s\ne\t'}} r_\s v_\s,\]
for some $r_\s\in\F$. Then
$\chi_\H^\lambda(T_w\tau)=\sum_\t\gamma_\t(w)$, so to determine the
character values $\chi_\H^\lambda(T_w\tau)$ it is enough to compute
$\gamma_\t(w)$, for $\t\in\Std(\lambda)$.

The key definition that we need comes from a paper of
Headley~\cite{Headley} who used it to compute the characters of the
alternating group $\Alt_n$. Two integers $i,j$, with $1\le i,j\le n$, are
\textbf{diagonally opposite} in~$\t$ if $\t(r,c)=i$ and $\t(c,r)=j$, for
some $(r,c)\in[\lambda]$ with $r\ne c$.  (Observe that if
$(r,c)\in\diag(\lambda)$ then $(c,r)\in\diag(\lambda)$ since
$\lambda=\lambda'$.)

\begin{Definition}[\protect{Headley~\cite[p.~130]{Headley}}]
Suppose that $w\in\Sym_n$ and that $\t$ is a standard $\lambda$-tableau,
for some partition $\lambda$. Then $\t$ is \textbf{$w$-transposable} if
whenever $i,j$ are diagonally opposite in~$\t$ then $i\in\{j\pm1\}$ and $i$ and
$j$ are in the same $w$--orbit.

Let $\Std(\lambda)_w$ be the set of $w$-transposable $\lambda$-tableaux.
\end{Definition}

Let $\Diag(\t)=\set{\t(r,r)|(r,r)\in\diag(\lambda)}$ be the set of numbers
that lie on the diagonal of the tableau~$\t$. If $\t$ is $w$-transposable and
$i\notin\Diag(\t)$ then $c_\t(i)=-c_\t(j)$ where $j\in\{i\pm1\}$ and $i$ and $j$
are diagonally opposite in $\t$. Consequently, the three cases in the next
definition are mutually exclusive (and exhaustive).

\begin{Definition}\label{D:gamma}
Suppose that $\lambda$ is a partition of~$n$ and that $w\in\Sym_n$ has a
reduced expression of the form $w=s_{i_1}\dots s_{i_k}$, where $1\le
i_1<\dots<i_k<n$. Let $\t$ be a $w$-transposable standard
$\lambda$-tableau and, for $j=1,\dots,k$, define
\[
  \gamma_\t(i_j)=\begin{cases}
      \frac{-1}{[\rho_j]}, &\text{if $i_j\in\diag(\t)$},\\
      \frac{-1}{[\rho'_j]}, &\text{if $c_\t(i_j)=-c_\t(i_j-1)$,}\\[5pt]
      \alpha_{\rho_j},&\text{if $c_\t(i_j)=-c_\t(i_j+1)$},\\
\end{cases}
\]
where $\rho_j=\rho_\t(i_j)=c_\t(i_j)-c_\t(i_j+1)$ and
$\rho'_j=c_\t(i_j-1)-c_\t(i_j+1)$.
\end{Definition}

Note that if $w\in\Sym_n$ then $w=s_{i_1}\dots s_{i_k}$, with $1\le
i_1<\dots<i_k<n$, if and only if $w=w_\kappa$ for some composition $\kappa$
of $n$. We have used reduced expressions in \cref{D:gamma}, rather
than compositions, only for convenience. See \cref{Example} below
for the definition in action.

\begin{Proposition}\label{factorization}
Suppose that $w=s_{i_1}\dots s_{i_k}$ is reduced, where
$1\le i_1<\dots<i_k<n$, and $\t\in\Std(\lambda)$. Then
\[
   \gamma_\t(w)=\begin{cases}
        \gamma_\t(i_1)\dots\gamma_\t(i_k),&
                \text{if $\t$ is $w$-transposable,}\\
        0,&\text{otherwise.}
\end{cases}
\]
\end{Proposition}

\begin{proof}
Observe that $\gamma_\t(w)\ne0$ (if and) only if there is a sequence
$\t_0,\t_1,\dots,\t_k$ of (not necessarily distinct) standard
$\lambda$-tableaux such that $\t_0=\t$, $\t_k=\t'$ and $v_{\t_j}$
appears with non--zero coefficient in $v_{\t_{j-1}}T_{i_j}$, for
$j=1,\dots,k$. By assumption, $1\le i_1<\dots<i_k<n$ so
$\gamma_\t(w)\ne0$ only if we can change $\t$ into $\t'$ by swapping
entries of the form $i_j$ and $i_j+1$, where $1\le j\le k$. Hence,
$\gamma_\t(w)\ne0$ only if $\t$ is $w$-transposable, proving
the first claim.

Suppose now that $\gamma_\t(w)\ne0$. Then the tableaux
$\t_0,\t_1,\dots,\t_k$ above are uniquely determined because $1\le i_1<\dots<i_k<n$.  Indeed, $\t_j$ is the unique standard
$\lambda$-tableau that has all of the numbers greater than $i_j+1$
in the same positions as they occur in~$\t$ and all numbers less than
or equal to $i_j+1$ in the same positions as in~$\t'$. It follows that
if~$v_{\t_j}$ appears with coefficient $\gamma_j$ in
$v_{\t_{j-1}}T_{i_j}$ then $\gamma_\t(w)=\gamma_1\dots\gamma_k$.
To complete the proof it remains to show that
$\gamma_\t(i_j)=\gamma_j$ is the coefficient of $v_{\t_j}$ in
$v_{\t_{j-1}}T_{i_j}$. There are three possibilities: either $i_j$ is
on the diagonal of $\t$, or $i_j$ is diagonally opposite either
$i_j-1$ or~$i_j+1$.

If $i_j\in\Diag(\t)$ then $\t_j=\t_{j-1}$ and, from the definitions,
$\gamma_j=\frac{-1}{[\rho_j]}=\gamma_\t(i_j)$, where
$\rho_j=\rho_{\t_j}(i_j)$.

Next, $i_j$ is diagonally opposite $i_j-1$ if and only if
$c_\t(i_j)=-c_\t(i_j-1)$. In this case we must have $i_{j-1}=i_j-1$
and $\t_{j-1}=\t_{j-2}s_{i_j-1}$ since otherwise we cannot swap
$i_j-1$ and $i_j$ in $\t$. Therefore,
$\gamma_j=\frac{-1}{[\rho_j']}=\gamma_\t(i_j)$, where
$\rho'_j=c_\t(i_j-1)-c_\t(i_j+1)$.

Finally, $i_j$ and $i_j+1$ are diagonally opposite in~$\t$ if and only
if $c_\t(i_j)=-c_\t(i_j+1)$. In this case we have
$\t_j=\t_{j-1}s_{i_j}$, so $\gamma_j=\alpha_{\rho_j}$, where
$\rho_j=c_\t(i_j)-c_\t(i_j+1)$. Hence, $\gamma_j=\gamma_\t(i_j)$ as
required.
\end{proof}

\begin{Corollary}\label{C:transposable}
  Suppose that $\lambda$ is a partition of $n$ and that
  $w=s_{i_1}\dots s_{i_k}$, reduced, with $1\le i_1\le\dots\le i_k<n$. Then
  \[\chi^\lambda(T_w\tau)
            =\sum_{\t\in\Std(\lambda)_w}\gamma_\t(i_1)\dots\gamma_\t(i_k).\]
  In particular, if there are no $w$-transposable $\lambda$-tableaux
  then $\chi^\lambda(T_w\tau)=0$.
\end{Corollary}

Hence, we need only consider $w_\kappa$-transposable tableaux to prove
\cref{T:TauCharacters}.

Recall from \cref{D:hlambda} that if $\lambda$ is a self-conjugate partition then
$h(\lambda)$ is the partition whose~$i^{\text{th}}$ part is the $(i,i)$--hook
length  of $\lambda$.  Our first aim is to show that
$\chi^\lambda(T_{w_\kappa}\tau)=0$ whenever $\kappa\ne h(\lambda)$.

\begin{Corollary}[\protect{cf. \cite[Lemma 3.3]{Headley}}]\label{C:odd}
Suppose that $\kappa$ is a composition of $n$ such that $\kappa$ has more
than $d(\lambda)$ parts of odd length.  Then
$\chi^\lambda(T_{w_\kappa}\tau)=0$.
\end{Corollary}

\begin{proof}
Observe that if $\t$ is a $w$-transposable tableau then each odd cycle
of $w$ meets the diagonal of $\t$ at least once. Consequently, there
are no transposable $w_\kappa$-tableaux because~$w_\kappa$ has more
than $d(\lambda)$ odd cycles. Hence, $\chi^\lambda(T_{w_\kappa}\tau)=0$
by \cref{C:transposable}.
\end{proof}

Before tackling other cycle types we prove a technical Lemma.

\begin{Lemma}\label{L:technical}
  Let $\kappa$ be a composition of $n$ and suppose that there exists a
  $w_\kappa$-transposable tableau $\t$ and integers $a$ and $b$ in
  $\Diag(\t)$ such $a<b$ are in the same cycle of $w_\kappa$
  and~$b$ is minimal with this property.  Let
  $\s=\t s_{a+1}s_{a+3}\dots s_{b-2}$.  Then $\s$ is
  $w_\kappa$-transposable and $\gamma_\s(w_\kappa)+\gamma_\t(w_\kappa)=0$.
\end{Lemma}

\begin{proof}Since $\t$ is $w_\kappa$-transposable, the numbers
  $\{a+2i+1,a+2i+2\}$, for $i=0,\dots,\frac{b-a-3}2$, are in the same
  $w_\kappa$-orbit and they occupy diagonally opposite positions in
  $\t$. Consequently, $b-a$ is odd and $\s$ is the tableau obtained from
  $\t$ by swapping these entries.  Hence, $\s$ is
  $w_\kappa$-transposable, proving our first claim.

Let $\gamma_\t[a,b)=\prod_{i=a}^{b-1}\gamma_\t(i)$ and define
$\gamma_{\s}[a,b)$ similarly. We claim that
$\gamma_\t[a,b)=-\gamma_{\s}[a,b)$. Establishing this claim will prove
that $\gamma_\t(w_\kappa)+\gamma_{\s}(w_\kappa)=0$ because
$\gamma_\t(i)=\gamma_{\s}(i)$ if $i<a$ or if $i\ge b$.

As in \cref{D:gamma}, set
$\rho_i=c_\t(i)-c_\t(i+1)=c_{\s}(i+1)-c_{\s}(i)$ and
$\rho'_i=c_\t(i-1)-c_\t(i+1)=c_{\s}(i+1)-c_{\s}(i-1)$, for
$a\le i<b$. Since $c_\t(a)=0=c_\s(a)$, if $i=0,\dots,b-a-1$ then
\begin{align*}
\gamma_\t(a+i)&=\begin{cases}
  \frac{-1}{[-c_\t(a+1)]},&\text{if $i=0$},\\
          \frac{-1}{[\rho_{a+i}']},&\text{if $i>0$ is even},\\
          \alpha_{\rho_{a+i}},&\text{if $i$ is odd},
\end{cases}
\intertext{and}
\gamma_{\s}(a+i)&=\begin{cases}
  \frac{-1}{[-c_\s(a+1)]},&\text{if $i=0$},\\
          \frac{-1}{[-\rho_{a+i}']},&\text{if $i>0$ is even},\\
          \alpha_{-\rho_{a+i}},&\text{if $i$ is odd}.
\end{cases}
\end{align*}
Let $c=\half(b-a-1)$ and recall from before \cref{D:SquareRoots} that
$[-k]=-q^{-2k}[k]$, for $k\in\Z$. Therefore,
\begin{align*}
  \gamma_{\t}[a,b)&=\frac{(-1)^{c+1}}{[-c_\t(a+1)]}
           \prod_{i=0}^{c-1}\frac{\alpha_{\rho_{a+2i+1}'}}{[\rho_{a+2i+2}']}
           =\frac{(-1)^{c+2}q^{2c_\t(a+1)}}{[c_\t(a+1)]}\prod_{i=0}^{c-1}
           \frac{q^{-2\rho'_{a+2i+2}}\alpha_{-\rho_{a+2i+1}'}}{[-\rho_{a+2i+2}']}\\
   &=-q^{2(c_\t(a+1)-\rho_{a+2}'-\dots-\rho_{b-1}')}\gamma_\s[a,b),
\end{align*}
since $c_\s(a+1)=-c_\t(a+1)$.  The definitions give a collapsing sum,
$c_\t(a+1)-\rho_{a+2}'-\dots-\rho_{b-1}'=c_\t(b)=0$. Hence,
$\gamma_\t[a,b)=-\gamma_{\s}[a,b)$. Therefore, as shown above,
$\gamma_\t(w_\kappa)+\gamma_{\s}(w_\kappa)=0$, proving the lemma.
\end{proof}

Recall that $\ell(\kappa)$ is the number of non-zero parts of~$\kappa$.

\begin{Corollary}[\protect{cf. \cite[Lemma 3.4]{Headley}}]\label{C:Many}
    Let $\kappa$ be a composition and suppose that $\ell(\kappa)<d(\lambda)$.
Then $\chi^\lambda(T_{w_\kappa}\tau)=0$.
\end{Corollary}

\begin{proof}
Let $\t$ be a $w_\kappa$-transposable tableau. As $w_\kappa$ has
$\ell(\kappa)<d(\lambda)$ cycles, at least one cycle in $w_\kappa$ meets the
diagonal of $\t$ more than once. Therefore, we can find integers $a<b$ on the
diagonal of $\t$, with $b$ minimal such that $a$ and $b$ are both in the same
cycle of $w_\kappa$. As in \cref{L:technical} let $\s=\t s_{a+1}s_{a+3}\dots
s_{b-2}$.  Then $\s$ is the $w_\kappa$-transposable obtained by swapping all of
the diagonally opposite pairs in~$\t$ that occur in the same cycle as~$a$
and~$b$. Therefore, the map $\t\mapsto\s$ gives a pairing of the
$w_\kappa$-transposable tableaux, so to prove the corollary it is enough to show
that $\gamma_\s(w_\kappa)+\gamma_{\t}(w_\kappa)=0$. However, this follows from
\cref{L:technical}.
\end{proof}

\begin{Corollary}[\protect{cf. \cite[Lemma 3.5]{Headley}}]\label{C:even}
Suppose that $\kappa$ is a composition of~$n$ and that $\kappa$ has at least
one even part. Then $\chi^\lambda(T_{w_\kappa}\tau)=0$.
\end{Corollary}

\begin{proof}Let $\t$ be a $w_\kappa$-transposable tableau and
consider the first even cycle in $\kappa$. As $\t$ is transposable,
this cycle meets the diagonal of $\t$ in an even number of places. If
this cycle meets the diagonal at $a<b$, where $b$ is minimal, then
$\gamma_\t(w_\kappa)+\gamma_{\s}(w_\kappa)=0$ by
\cref{L:technical}, where $\s=\t s_{a+1}\dots s_{b-2}$. If the
first even cycle $(b-1,b-2,...,a+1)$ does not meet the diagonal of
$\t$ at all then we again define $\s=\t s_{a+1}\dots s_{b-2}$. By
essentially repeating the argument of \cref{L:technical} we find
that $\gamma_\t(w_\kappa)+\gamma_{\s}(w_\kappa)=0$. Hence, it
follows that $\chi^\lambda(T_{w_\kappa}\tau)=0$ as required.
\end{proof}

To prove \cref{T:TauCharacters} it remains to consider
those elements $w_\kappa$ that have $d(\lambda)$
cycles of odd length, and no even cycles. This is by far the most
complicated case and it requires some new combinatorial
machinery. We begin by reformulating an identity of Greene's~\cite{Green}.

Let $(X,<)$ be a poset. Then $X$ is \textbf{connected} if its Hasse
diagram is connected; otherwise $X$ is \textbf{disconnected}. If
$X$ contains $m+1$ elements then a
\textbf{linearisation} of $X$ is a bijection $f\map X\{0,\dots,m\}$ that
respects the ordering in~$X$; thus, $f(x)<f(y)$ whenever $x<y$ for
$x,y\in X$.  Let $\L(X)$ be the set of linearisations of $X$. If
$f\in\L(X)$ let $f^{-1}\map{\{0,\dots,m\}}X$ be its inverse.
A poset $(X,<)$ is \textbf{semilinear} if, in the Hasse diagram of~$X$,
every vertex has at most two edges. Write $x\lessdot y$ if $x<y$ and $x$
and $y$ are adjacent in $X$ (that is, $x\le a\le y$ only if $a=x$ or $a=y$).

Let $X$ be a semilinear poset with $m+1$ elements. Fix a labelling
$X=\{x_0,x_1,\dots,x_m\}$ of the elements of $X$ so that $x_i\lessdot x_j$
only if $j\in\{i\pm1\}$. If $f\in\L(X)$ write $f^*(i)=j$ if
$f^{-1}(i)=x_j$, for $0\le i\le m$.  Define
\begin{equation}\label{E:PosetSign}
\epsilon_X(i)=\begin{cases}
       1,& \text{if } x_i\lessdot x_{i+1},\\
       -1,&\text{if } x_{i+1}\lessdot x_i,\\
       0,&\text{otherwise},
\end{cases}
\end{equation}
for $0\le i<m$, and set $\epsilon_X=\prod_{i=0}^{m-1}\epsilon_X(i)$. Then~$X$ is
disconnected if and only if $\epsilon_X=0$.

The following result is a mild reformulation of
Greene~\cite[Theorem~3.4]{Greene}.

\begin{Lemma}[Greene~\cite{Greene}]\label{L:Greene}
Suppose that $(X,<)$ is a semilinear poset with elements $\{x_0,\dots,x_m\}$
labelled as above.  Suppose that $\{c_0,\dots,c_m\}$ is a set of pairwise
distinct integers. Then
\[
\sum_{f\in\L(X)}q^{2c_{f^*(m)}}\prod_{i=0}^{m-1} \frac{1}{[c_{f^*(i+1)}-c_{f^*(i)}]}
=\epsilon_Xq^{2c_m}\prod_{i=0}^{m-1} \frac{1}{[c_{i+1}-c_i]}.
\]
\end{Lemma}

\begin{proof}Let $\set{Q_x|x\in X}=\{Q_0,\dots,Q_m\}$ be a set of
  indeterminates.  Greene~\cite[Theorem~3.4]{Greene} has shown that
  \[ \sum_{f\in\L(X)}\prod_{i=0}^{m-1} \frac1{Q_{f^*(i+1)}-Q_{f^*(i)}}
           =\epsilon_X\prod_{i=0}^{m-1}\frac1{Q_{i+1}-Q_i}.
  \]
  To deduce the Lemma from Greene's identity set $Q_x=q^{2c_x}$, multiply
  both sides by $q^{2(c_0+\dots+c_m)}(q-q^{-1})^m$ and then simplify each
  of the factors using the observation that
  $\frac{q^{2b}(q-q^{-1})}{q^{2a}-q^{2b}}=\frac1{[a-b]}$, for $a,b\in\Z$.
\end{proof}

\begin{Remark}
  In fact, Greene~\cite{Greene} proves a more general product formula
  for \textit{planar} posets. When comparing \cref{L:Greene} with
  Greene's theorem note that because $X$ is semilinear if $x\ne y\in X$
  then the value of M\"obius function~$\mu$ at~$(x,y)$ is given by
  \[  \mu(x,y)=\begin{cases}
              -1, &\text{if }x\lessdot y,\\
              0,  &\text{otherwise}.
      \end{cases}
   \]
   The scalar $\epsilon_X$ does not appear in \cite{Greene}. We
   introduce $\epsilon_X$ as a device for keeping track of the order of
   the terms on the left-hand side of the identity in \cref{L:Greene}.
   These signs are a necessary source of pain below.
\end{Remark}

% If $\lambda$ and $\nu$ are partitions with
% $\Diag(\nu)\subset\Diag(\lambda)$ then a \textbf{skew}
% $\lambda/\nu$-tableau is an injective map
% $\s\map{\Diag(\lambda)/\Diag(\nu)}\Z$. If $\t$ is a
% $\lambda$-tableau then a \textbf{subtableau} of $\t$ is any map
% $\s\map D\Z$, where $D\subseteq\diag(\lambda)$. For $i=1,\dots,d(\lambda)$
% let $\t_i$ be the subtableau of $\t$ that contains the numbers in the
% $i^\th$ diagonal hook $H_i(\lambda)$; that is, $\t_i$ is the restriction
% of~$\t$ to $H_i(\lambda)$.  We think of~$\t_i$ as being a tableau of shape
% $(\frac{h_i+1}2,1^{\frac{h_i-1}2})$ in the natural way.

We now engineer the circumstances that we need to apply \cref{L:Greene}.

Until further notice, fix a composition
$\kappa=(\kappa_1,\dots,\kappa_d)$, such that
$\kappa_1,\dots,\kappa_d$ are all odd and $d=d(\lambda)$.
Let~$\t$ be a $w_\kappa$-transposable tableau~$\t$.
Then $\t$ determines a sequence of partitions
$\lambda_{\t,0}=(0),\lambda_{\t,1},\dots,\lambda_{\t,d}=\lambda$ where
$(r,c)\in\diag(\lambda_{\t,z})$ if and only if
$\t(r,c)\le\kappa_1+\dots+\kappa_z$, for $1\le z\le d$. That is,
$\lambda_{\t,z}{/}\lambda_{\t,z-1}$ is the shape of the skew subtableau of
$\t$ that contains the numbers in cycle $z$ of $w_\kappa$. Each of the
partitions $\lambda_{\t,1},\dots,\lambda_{\t,z}$ is symmetric so,
in particular, $(z,z)\in\diag(\lambda_{\t,z})$, for all $z$.

Fix an integer $z$, with $1\le z\le\ell(\kappa)=d$, and let
$X_{\t,z}=\set{(r,c)\in\diag(\lambda_{\t,z}/\lambda_{\t,z-1})|c\ge r}$,
which we consider as a partially ordered set with $(r,c)\preceq (r',c')$
if $r\le r'$ and $c\le c'$. Set $k_z=\kappa_1+\dots+\kappa_{z-1}+1$ and
$m_z+1=|X_{\t,z}|=\half(\kappa_z-1)$, so that $(k_z,k_z+1,\dots,k_z+2m_z)$ is the $z$th
cycle in $w_\kappa$. If $x=(r,c)\in X_{\t,z}$ then set
$x'=(c,r)$. Consequently, if $x\in X_{\t,z}$ then, $x'\in X_{\t,z}$ if and only
if $x=(z,z)$. Set $\omega=(z,z)$ and define two closely related \textit{sign sequences}
$\epsilon_{\t,z}\map{X_{\t,z}}\{\pm1\}$ and
$\epsilon^\omega_{\t,z}\map{X_{\t,z}}\{\pm1\}$ by
\[\epsilon_{\t,z}(x)=\begin{cases} 1,&\text{if }\t(x')\ge\t(x),\\ -1,&\text{if
  }\t(x')<\t(x), \end{cases} \quad\text{and}\quad
  \epsilon^\omega_{\t,z}(x)=\begin{cases} 1,&\text{if }\t(x)\ge\t(\omega),\\ -1,&\text{if
    }\t(x)<\t(\omega), \end{cases}
\]
for $x\in X_{\t,z}$. By definition,
$\epsilon_{\t,z}(\omega)=1=\epsilon^\omega_{\t,z}(\omega)$ so we could, instead,
define the sign sequences to be functions from
$X_{\t,z}\setminus\set{\omega}\longrightarrow\set{\pm1}$. When $z$ is
understood we omit it and simply write $\epsilon_\t$ and
$\epsilon^\omega_\t$. For the tableau $\t$, define a linearisation
$f_{\t,z}\in\L(X_{\t,z})$ of $X_{\t,z}$ by
\[ f_{\t,z}(x)=\#\set{y\in X_{\t,z}|\t(y)\le\t(x)}, \]
%\[f_{\t,z}(x)=\#\set{x\in X_{\t,z}|\t(x)\ge\t(\omega)}
%         =\#\set{x\in X_{\t,z}|\epsilon^\omega_{\t,z}(x)=1},\]
for $x\in X_{\t,z}$.

As a prelude to applying \cref{L:Greene}, the next definition partitions
$\Std(\lambda)_{w_\kappa}$ into more manageable pieces.

\begin{Definition}\label{D:SimEquiv}
  Suppose that $\kappa$ is a composition of~$n$ and that $1\le z\le d(\lambda)$.
  Let~$\sim_z$ be the equivalence relation on the set of $w_\kappa$-transposable
  $\lambda$-tableau determined by $\s\sim_z\t$, for
  $\s,\t\in\Std(\lambda)_{w_\kappa}$, if:
\begin{enumerate}
  \item $X_{\s,z}=X_{\t,z}$ (so that $\lambda_{\s,z}=\lambda_{\t,z}$),
  \item $\s(x)=\t(x)$, for all $x\in\diag(\lambda)/X_{\s,z}$, and,
  \item $\epsilon_{\s,z}(x)\epsilon^\omega_{\s,z}(x)
  =\epsilon_{\t,z}(x)\epsilon^\omega_{\t,z}(x)$, for all $x\in X_{\s,z}$.
\end{enumerate}
If $\T$ is a $\sim_z$-equivalence class of $w_\kappa$-transposable tableaux set
$X_\T=X_{\t,z}$, $\lambda_\T=\lambda_{\t,z}$ and
\[
  \varepsilon_\T=\prod_{x\in X_\T}\varepsilon_\T(x), \qquad\text{where}\quad
  \varepsilon_\T(x)=\epsilon_{\t,z}(x)\epsilon_{\t,z}^\omega(x),
\]
and $\t$ is any element of $\T$.
\end{Definition}

The definition of the $\sim_z$-equivalence relation ensures that $X_\T$,
$\lambda_\T$ and $\varepsilon_\T$ depend only on~$\T$ and not on the
choice of~$\t\in\T$. We warn the reader that, in general, the signs
$\epsilon_{X_\T}$ and $\varepsilon_\T$ can be different.
% Note that
% $\varepsilon_\T(\omega)=\epsilon_\t(\omega)\epsilon^\omega_\t(\omega)=1$,
% for all $\t\in\Std(\lambda)_{w_\kappa}$, so $x=\omega$ can be omitted
% from the definition of $\varepsilon_\T$.

As the next result suggests, the $\sim_z$-equivalence classes will allow us to
apply \cref{L:Greene}.

\begin{Lemma}\label{L:Linearisation}
  Suppose that $1\le z\le d(\lambda)$ and let $\T$ be a
  $\sim_z$-equivalence class of $w_\kappa$-transposable $\lambda$-tableaux.
  Then the map $\t\mapsto f_{\t,z}$ defines a bijection $\T\bijection\L(X_\T)$.
\end{Lemma}

\begin{proof}
  To show that the map $\T\mapsto\L(X_\T)$ is a bijection we define an
  inverse map. First note that if $\s,\t\in\T$ then $\t$ and~$\s$ agree
  on~$\lambda\setminus\lambda_\T$. Therefore, if~$f\in\L(X_\T)$ then to define a
  tableau $\t_f\in\T$ we only need to specify its values on
  $\lambda/\lambda_\T$. All but one of the numbers
  $k_z,k_z+1,\dots,k_z+2m_z$ occurs in a diagonally opposite pair, so
  the requirement that
  $\varepsilon_\T(x)=\epsilon_{\t_f}(x)\epsilon^\omega_{\t_f}(x)$ is
  constant on~$\T$ implies that there is a unique tableau $\t_f\in\T$
  such that $\t_f(\omega)=k_z+2f(\omega)$ and if $x\in X_\T$ and
  $x\ne\omega=(z,z)$ then
  \[\{\t_f(x),\t_f(x')\}=\begin{cases}
        \{k_z+2f(x),k_z+2f(x)+1\}, &\text{if }f(x)<f(\omega),\\
        \{k_z+2f(x)-1,k_z+2f(x)\}, &\text{if }f(x)>f(\omega).
      \end{cases}
  \]
  The maps $\t\mapsto f_\t$ and $f\mapsto\t_f$ are mutually inverse, so
  $\T\bijection\L(X_\T)$ as claimed.
\end{proof}

The following example should help the reader absorb all of the new definitions.

\begin{Example}\label{Example} Suppose that $\lambda=(6,3,2,1^3)$, so
  that $h(\lambda)=(11,3)$ and take $\kappa=(7^2)$. Then
  $w_\kappa=s_1\dots s_6s_8\dots s_{13}$ and the reader may check that
  there are $384$ $w_\kappa$-transposable $\lambda$-tableaux.  We take
  $z=2$, so that $k_z=8$ and $m_z=3$, and consider the
  $\sim_z$-equivalence class~$\T$ of $w_\kappa$-transposable tableaux
  with $\lambda_{\t,z}=\lambda/(4,1^3)$ and with sign sequence
  $\varepsilon_\T={-}{-}{+}$. That is,
  $\epsilon_{\t,z}(x_1)\epsilon^\omega_{\t,z}(x_2)=-1$,
  $\epsilon_{\t,z}(x_2)\epsilon^\omega_{\t,z}(x_3)=-1$ and
  $\epsilon_{\t,z}(x_3)\epsilon^\omega_{\t,z}(d)=1$, where $x_0=(2,2)=\omega$,
  $x_1=(2,3)$, $x_2=(1,5)$ and $x_3=(1,6)$. We use this shorthand for all sign
  sequences in the table below.
  \[\begin{array}{*3c@{\ }*7{@{\ }c}}
  \t&\epsilon_{\t,2}&\epsilon^\omega_{\t,2}&f_{\t,2}^*
     &\gamma_\t(8)&\gamma_\t(9)&\gamma_\t(10)&\gamma_\t(11)
     &\gamma_\t(12)&\gamma_\t(13)\\\toprule
  \Tableau{{\cdot,\cdot,\cdot,\cdot,12,13},{\cdot,8,10},{\cdot,9},\cdot,11,14}
     &{-}{-}{+}&{+}{+}{+}&0123
     &\frac{-1}{[0+1]}&\alpha_{2}&\frac{-1}{[-1+4]}&\alpha_{8}
     &\frac{-1}{[-4-5]}&-\alpha_{10}\\[9mm]
  \Tableau{{\cdot,\cdot,\cdot,\cdot,10,13},{\cdot,8,12},{\cdot,11},\cdot,9,14}
     &{-}{-}{+}&{+}{+}{+}&0213
     &\frac{-1}{[0+4]}&\alpha_{8}&\frac{-1}{[-4+1]}&\alpha_{2}
     &\frac{-1}{[-1-5]}&-\alpha_{10}\\[9mm]
  \Tableau{{\cdot,\cdot,\cdot,\cdot,10,11},{\cdot,8,14},{\cdot,13},\cdot,9,12}
     &{-}{-}{+}&{+}{+}{+}&0312
     &\frac{-1}{[0+4]}&\alpha_{8}&\frac{-1}{[-4-5]}&-\alpha_{10}
     &\frac{-1}{[5+1]}&\alpha_{2}\\[9mm]
  %\end{array}\]
  %\[\begin{array}{*3c@{\ }*7{@{\ }c}}
  \Tableau{{\cdot,\cdot,\cdot,\cdot,8,13},{\cdot,10,12},{\cdot,11},\cdot,9,14}
     &{-}{+}{+}&{+}{-}{+}&1203
     &-\alpha_{8}&\frac{-1}{[4-0]}&\frac{-1}{[0+1]}&\alpha_{2}
     &\frac{-1}{[-1-5]}&-\alpha_{10}\\[9mm]
  \Tableau{{\cdot,\cdot,\cdot,\cdot,8,11},{\cdot,10,14},{\cdot,13},\cdot,9,12}
     &{-}{+}{+}&{+}{-}{+}&1302
     &-\alpha_{8}&\frac{-1}{[4-0]}&\frac{-1}{[0-5]}&-\alpha_{10}
     &\frac{-1}{[5+1]}&\alpha_{2}\\[9mm]
  \Tableau{{\cdot,\cdot,\cdot,\cdot,8,11},{\cdot,12,14},{\cdot,13},\cdot,9,10}
     &{-}{+}{-}&{+}{-}{-}&2301
     &-\alpha_{8}&\frac{-1}{[4+5]}&\alpha_{10}&\frac{-1}{[-5-0]}
     &\frac{-1}{[0+1]}&\alpha_{2}\\[9mm]
  \bottomrule
  \end{array}\]
  In the tableaux above we have only indicated the positions of the numbers $8,\dots,14$
  because these are the only entries that will matter in the arguments below
  (because the tableaux in a $\sim_z$-equivalence class are constant on the
  nodes outside of $X_\T$).
\end{Example}

Now fix a $\sim_z$-equivalence class of tableaux~$\T$ with
poset $X_\T$, where $1\le z\le d(\lambda)$. Define
\begin{equation}\label{E:alphaT}
  \alpha(X_\T)=\prod_{x\in X_\T\setminus\set{\omega}} \alpha_{2c(x)}.
\end{equation}
For a $w_\kappa$-transposable tableau $\t$ set
\begin{equation}\label{E:gammatz}
\gamma_{\t,z}(w_\kappa)
    =\prod_{j=0}^{2m_z-1}\gamma_\t(k_z+j).
\end{equation}
Then, $\gamma_\t(w_\kappa)=\gamma_{\t,1}(w_\kappa)\dots\gamma_{\t,d}(w_\kappa)$.
The reader may find it helpful to consult \cref{Example} above during the
proofs of the next few results.

\begin{Lemma}\label{L:GammaReduction}
    Suppose that $\kappa=(\kappa_1,\dots,\kappa_d)$ is a composition of
    $n$ such that $d=d(\lambda)$ and $\kappa_1,\dots,\kappa_d$ are all
    odd.  Fix $1\le z\le d$ and let $\T$ be a $\sim_z$-equivalence class
    of $w_\kappa$-transposable tableaux and $\t\in\T$. Then
    \[
        \gamma_{\t,z}(w_\kappa)=\dfrac{\alpha(X_\T)\Prod_{i=0}^{m_z}\epsilon_\t(x^\t_i)}
                                       {\Prod_{i=0}^{m_z-1}[c_\t(i)-c_\t(i+1)]},
    \]
    where $x^\t_i=f_{\t,z}^{-1}(i)$ and $c_\t(i)=\epsilon_\t(x^\t_i)c(x^\t_i)$,
    for $i=0,\dots,m_z$.
\end{Lemma}

\begin{proof}
  By definition, $\gamma_{\t,z}(w_\kappa)=\prod_{j=0}^{2m_z-1}\gamma_\t(k_z+j)$
  and $\alpha(X_\T)=\prod_{x\in X_\T} \alpha_{2c(x)}$. So,
  we need to show that
  \[ \prod_{j=0}^{2m_z-1}\gamma_\t(k_z+j)
        = \prod_{x\in X_\T} \epsilon_\t(x)\alpha_{2c(x)}\prod_{i=0}^{m_z-1}\cdot
             \frac1{[\epsilon_\t(x^\t_i)c(x^\t_i)-\epsilon_\t(x^\t_{i+1})c(x^\t_{i+1})]}.
  \]
  There are three mutually exclusive cases to consider, corresponding
  to the different cases in \cref{D:gamma}.

  Fix~$j$, with $1\le j\le 2m_z$, and let $x$ and $y$ be the unique nodes in
  $X_\t$ such that $k_z+j\in\{\t(x),\t(x')\}$ and
  $k_z+j+1\in\{\t(y),\t(y')\}$, where we allow and, in fact, need to
  include the possibility that $\set{x,x'}=\set{y,y'}$.

  \textit{Case 1.} $k_z+j\in\Diag(\t)$. Therefore, $x=(z,z)=\omega$ and
  $c_\t(k_z+j)=0$. Write $x=x^\t_i$, so that $y=x^\t_{i+1}$. Now, $\t(y)\ge\t(x)+1$ if
  and only if $\epsilon_\t(y)=1$, in which case $c_\t(k_z+j+1)=c(x^\t_{i+1})$.
  Similarly, $\t(y')\ge\t(x)+1$ if and only if $\epsilon_\t(y)=-1$ in which
  case $c_\t(k_z+j+1)=-c(x^\t_{i+1})$. Hence, in both cases,
  $c_\t(k_z+j+1)=\epsilon_\t(x^\t_{i+1})c(x^\t_{i+1})$.  Therefore,
  \[\gamma_\t(k_z+j)=\frac{-1}{[c_\t(k_z+j)-c_\t(k_z+j+1)]}
  =\frac{-1}{[\epsilon_\t(x^\t_i)c(x^\t_i)-\epsilon_\t(x^\t_{i+1})c(x^\t_{i+1})]}.\]

  \textit{Case 2.} $c_\t(k_z+j)=-c_\t(k_z+j-1)$. In this case $k_z+j$
  and $k_z+j-1$ are diagonally opposite. Hence, if
  $x=x^\t_{i+1}$ then $y=x^\t_i$.
  Further, $\t(x)=k_z+j-1$ if and only if $\epsilon_\t(x)=1$, and
  $\t(y)=k_z+j-1$ if and only if $\epsilon_\t(y)=1$. Hence, as in Case~1,
  \[\gamma_\t(k_z+j)
      = \frac{-1}{[c_\t(k_z+j-1)-c_\t(k_z+j+1)]}
      =\frac{-1}{[\epsilon_\t(x^\t_i)c(x^\t_i)-\epsilon_\t(x^\t_{i+1})c(x^\t_{i+1})]}.
  \]

  \textit{Case 3.} $c_\t(k_z+j)=-c_\t(k_z+j+1)$. In this case $k_z+j$
  and $k_z+j+1$ are diagonally opposite so~$x=y$ and, in particular,
  $x\ne (z,z)$. Write $x=x^\t_i$. Then
  \[\gamma_\t(k_z+j)
      =\alpha_{c_\t(k_z+j+1)-c_\t(k_z+j)}
      =\alpha_{-2\epsilon_\t(x^\t_i)c(x^\t_i)}
      =-\epsilon_\t(x^\t_i)\alpha_{2c(x^\t_i)},\]
  where the last equality uses the fact that $\alpha_{-c}=-\alpha_c$, for $c\in\Z$.

  By construction, $0\le i<m_z$ in Cases~1--3 with
  $\omega\in\set{x^\t_i,x^\t_{i+1}}$ only in Case~1. Cases~1 and~2
  contribute the same quantity, as a function of~$i$, to $\gamma_{\t,z}(w_\kappa)$.
  There are $2m_z$ minus signs in Cases~1--3, so
  the lemma now follows by combining these three cases since
  $\prod_{x\ne\omega}\epsilon_\t(x)=\prod_{x\in X_\T}\epsilon_\t(x)$.
\end{proof}

As before \cref{L:Greene}, fix a labelling $X_\T=\{x_0,x_1,\dots,x_{m_z}\}$ that
is compatible with the partial order on~$X_\T$ in the sense that
$x_i\lessdot x_j$ only if $j\in\{i\pm1\}$. For $0\le i\le m_z$ define
\[ c_\T(i)=\varepsilon_\T(x_i)c(x_i)\in\Z. \]
We are now ready to compute $\sum_{\t\in\T}\gamma_{\t,z}(w_\kappa)$.

\begin{Proposition}\label{P:transposableReduction}
  Suppose that $\kappa=(\kappa_1,\dots,\kappa_d)$ is a composition of $n$
  that has $d=d(\lambda)$ odd parts. Fix $z$ with
  $1\le z\le d(\lambda)$ and let $\T$ be a $\sim_z$-equivalence class of
  $w_\kappa$-transposable $\lambda$-tableaux. Then
  \[\sum_{\t\in\T}\gamma_{\t,z}(w_\kappa)=\varepsilon_\T\epsilon_{X_\T}(-1)^{m_z}
       q^{2c_\T(m_z)}\alpha(X_\T)\, \prod_{i=0}^{m_z-1}\frac1{[c_\T(i+1)-c_\T(i)]}.
  \]
\end{Proposition}

\begin{proof}
As in \cref{L:GammaReduction}, for $\t\in\T$ let
$c_\t(i)=\epsilon_\t(x_i^\t)c(x_i^\t)$, where $x_i^\t=f_\t^{-1}(i)$ and
$0\le i\le m_z$. Then
\[\sum_{\t\in\T}\gamma_\t(w_\kappa)
=\alpha(X_\T)\sum_{\t\in\T}\frac{\prod_{i=0}^{m_z}\epsilon_\t(x_i)}
                 {\prod_{i=0}^{m_z-1}[c_\t(i)-c_\t(i+1)]}.
\]
by \cref{L:GammaReduction}. Now, if $\t\in\T$ then $x_i^\t=x_j$, where
$j=f^*_\t(i)$. So, using \cref{D:SimEquiv}(c),
\[c_\t(i)=\epsilon_\t(x_i^\t) c(x_i^\t)
   =\epsilon^\omega_\t(x_j)\epsilon^\omega_\t(x_j)\epsilon_\t(x_j)c(x_j)
   =\epsilon^\omega_\t(x_j)c^\T_j
   =\epsilon^\omega_\t(x^\t_i)c_\T(f_\t^*(i)).
\]
By definition, $\epsilon^\omega_\t(x^\t_i)=-1$ if and only if
$0\le i<f_\t(\omega)$.  Moreover,
$\epsilon_\t(x)=\varepsilon_\T(x)\epsilon^\omega_\t(x)$, for $x\in X_\T$.
Therefore, using the fact that $c(\omega)=0$ for the second and fourth
equalities,
\begin{align*}
  \sum_{\t\in\T}\gamma_\t(w_\kappa)
    &=\displaystyle\alpha(X_\T)\sum_{\t\in\T}
      \frac{\prod_{i=0}^{m_z}\varepsilon_\T(x_i)\epsilon^\omega_\t(x_i)}%
           {\prod_{i=0}^{m_z-1}[\epsilon^\omega_\t(x^\t_i)c_\T(f_\t^*(i))
               -\epsilon^\omega_\t(x^\t_{i+1})c_\T(f_\t^*(i+1))]}\\
    &=\displaystyle\varepsilon_\T\alpha(X_\T)\sum_{\t\in\T}
    \prod_{i=0}^{f_\t(\omega)-1}\frac{-1}{[c_\T(f_\t^*(i+1))-c_\T(f_\t^*(i))]}
    \prod_{i=f_\t(\omega)}^{m_z-1}\frac{1}{[c_\T(f_\t^*(i))-c_\T(f_\t^*(i+1))]}\\
    &=\displaystyle\varepsilon_\T\alpha(X_\T)\sum_{\t\in\T}
    \prod_{i=0}^{f_\t(\omega)-1}\frac{-1}{[c_\T(f_\t^*(i+1))-c_\T(f_\t^*(i))]}
    \prod_{i=f_\t(\omega)}^{m_z-1}\frac{-q^{2(c_\T(f_\t^*(i+1))-c_\T(f_\t^*(i)))}}%
                {[c_\T(f_\t^*(i+1))-c_\T(f_\t^*(i))]}\\
    &=\displaystyle\varepsilon_\T(-1)^{m_z}\alpha(X_\T)
       \sum_{\t\in\T}q^{2c_\T(f_\t^*(m_z))}
         \prod_{i=0}^{m_z-1}\frac{1}{[c_\T(f_\t^*(i+1))-c_\T(f_\t^*(i))]}\\
    &=\displaystyle\varepsilon_\T(-1)^{m_z}\alpha(X_\T)
      \sum_{f\in\L(X_\T)}q^{2c_\T(f^*(m_z))}
         \prod_{i=0}^{m_z-1}\frac{1}{[c_\T(f^*(i+1))-c_\T(f^*(i))]},
\end{align*}
where the last equation invokes the bijection $\T\bijection \L(X_\T)$ of
\cref{L:Linearisation}. Applying \cref{L:Greene} now completes the proof.
\end{proof}

If the poset $X_\T$ is disconnected then $\epsilon_{X_\T}=0$ so that
$\sum_{\t\in\T}\gamma_{\t,z}(w_\kappa)=$ by
\cref{P:transposableReduction}. We therefore need to determine when
$X_\T$ is connected.

Let $\lambda$ and $\mu$ be partitions with $\mu\subset\lambda$. Then
$\lambda/\mu$ is a \textbf{strip} if the multiset
$\set{c(x)|x\in\diag(\lambda/\mu)}$ contains distinct consecutive
integers.

\begin{Definition}
Suppose that $\kappa=(\kappa_1,\dots,\kappa_d)$ is a
composition of~$n$. A \textbf{symmetric covering} of $\lambda$ of
\textbf{type} $\kappa$ is a sequence of self-conjugate partitions
\[ \lambda^{(0)}=(0)\subset\lambda^{(1)}\subset\dots\subset\lambda^{(d)}=\lambda \]
such that $\lambda^{(z)}/\lambda^{(z-1)}$ is a strip and
$|\lambda^{(z)}|=\kappa_1+\dots+\kappa_z$, for $1\le z\le d$.
The composition~$\kappa$ \textbf{symmetrically covers} $\lambda$ if a
symmetric covering of $\lambda$ of type $\kappa$ exists.
\end{Definition}

For example, if
$\lambda=(4,3,3,1)$ and $\kappa=(3,1,7)$ then the following diagrams show
that~$\kappa$ symmetrically covers $\lambda$:
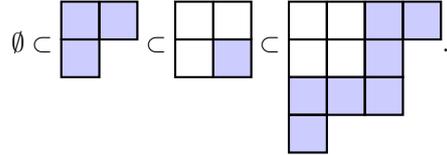
\[\emptyset\subset\ShadedTableau[(1,0),(2,0),(1,-1)]{{,},{\space}}
           \subset\ShadedTableau[(2,-1)]{{,},{, }}
           \subset\ShadedTableau[(1,-3),(1,-2),(2,-2),(3,-2),(3,-1),(3,0),(4,0)]%
                   {{,,,},{,,},{,,},\ }.
\]
Notice that $h(\lambda)=(7,3,1)$ so that $\vec\kappa=h(\lambda)$.

\begin{Lemma}\label{L:covers}
  Suppose that $\kappa$ is a composition of $n$ such that
  $\len(\kappa)=d(\lambda)$. Then $\kappa$ symmetrically
  covers $\lambda$ if and only if $\vec\kappa= h(\lambda)$. Moreover, if
  $\vec\kappa=h(\lambda)$ then there is a unique symmetric covering of
  $\lambda$ of type $\kappa$.
\end{Lemma}

\begin{proof}
Write $\kappa=(\kappa_1,\dots,\kappa_d)$, where $d=d(\lambda)=\len(\kappa)$.
We argue by induction on $d$. If $d=1$ then $\kappa$ symmetrically covers
$\lambda$ if and only if $n$ is odd, $\kappa=(n)=h(\lambda)$ and
$\lambda=(\frac{n+1}2,1^{\frac{n-1}2})$ is a hook partition. Now suppose that
$d>1$ and write $h(\lambda)=(h_1,\dots,h_d)$.

Suppose that $\kappa$ symmetrically covers $\lambda$. Then
$\nu=(\kappa_1,\dots,\kappa_{d-1})$ symmetrically covers $\lambda^{(d-1)}$ so
that $\vec\nu=h(\lambda^{(d-1)})$ by induction. As $\lambda=\lambda^{(d)}$ and
$\lambda^{(d-1)}$ are both symmetric and $\lambda^{(d)}/\lambda^{(d-1)}$ is
connected, it follows that $\kappa_d$ is odd and that $\kappa_d=h_i$ where $i$
is minimal such that $\lambda_i\ne\lambda^{(d-1)}_i$. Hence, $\vec\kappa=
h(\lambda)$.

Conversely, suppose that $d>1$ and that $\vec\kappa= h(\lambda)$ and let
$\nu=(\kappa_1,\dots,\kappa_{d-1})$. Then $\kappa_d=h_i$ for some $i$, so
that $\vec{\nu}=h(\mu)=(h_1,\dots,h_{i-1},h_{i+1},\dots,h_d)$, where the
diagram of $\mu$ is obtained by deleting the $(i,i)$-rim hook from the
diagram of $\lambda$ (the $(i,i)$-rim hook is the set of nodes along the
rim of $\lambda$ that connect the nodes $(i,\lambda_i)$ and
$(\lambda_i,i)$).

Finally, the uniqueness of the covering is immediate from the
construction in the last paragraph.
\end{proof}

Combining \cref{L:covers} with \cref{P:transposableReduction} we can
prove that the characters $\chi^\lambda(T_{w_\kappa}\tau)$ vanish when
$\kappa$ has $d(\lambda)$ odd parts and $\vec\kappa\ne h(\lambda)$.

\begin{Corollary}\label{C:zero!}
Let $\kappa=(\kappa_1,\dots,\kappa_d)$ be a composition of $n$
that has $d=d(\lambda)$ non-zero parts, all of which are odd, and suppose
that $\vec\kappa\ne h(\lambda)$. Then
$\chi^\lambda_\H(T_{w_\kappa}\tau)=0$.
\end{Corollary}

\begin{proof} By \cref{C:transposable} it is enough to show that
  $\sum_\t\gamma_\t(w_\kappa)=0$, where $\t$ runs over the set of
  $w_\kappa$-transposable $\lambda$-tableaux.
  Let $\t$ be a $w_\kappa$-transposable $\lambda$-tableau. By assumption
  $\vec\kappa\ne h(\lambda)$ so, by \cref{L:covers}, we can find an
  integer~$z$ such that $1\le z\le d(\lambda)$ and
  $\lambda_{\t,z}/\lambda_{\t,z-1}$ is not connected. Let~$z$ be the
  smallest such integer and let  $\T$ be the $\sim_z$-equivalence class
  of $w_\kappa$-transposable $\lambda$-tableaux that contains~$\t$. In
  this way we attach a pair $(z,\T)$ to each $w_\kappa$-transposable
  $\lambda$-tableau (different $z$'s can appear for different tableaux).
  The equivalence classes that we obtain in this way partition the set
  of $w_\kappa$-transposable $\lambda$-tableaux because the choices of
  $\sim_z$-equivalence classes are determined by the cycles of
  $w_\kappa$ and the minimality of~$z$. Hence, it is enough to show that
  $\sum_{\t\in\T}\gamma_\t(w_\kappa)=0$ whenever $\T$ is one of the
  chosen $\sim_z$-equivalence classes of $w_\kappa$-transposable
  $\lambda$-tableau.  The poset $X_\T$ is disconnected by the choice of
  $z$.  Therefore, $\epsilon_{X_\T}=0$ and so
  $\sum_{\t\in\T}\gamma_\t(w_\kappa)=0$ by
  \cref{P:transposableReduction}.
\end{proof}

We have now computed $\chi^\lambda(T_{w_\kappa\tau})$ in all cases except for
when $\vec\kappa=h(\lambda)$. Fix a composition $\kappa$ such
that $\vec\kappa=h(\lambda)$ together with the symmetric covering
$\lambda^{(0)}\subset\dots\subset\lambda^{(d)}$ of~$\lambda$ of type~$\kappa$
given by \cref{L:covers}.

\begin{Corollary}\label{C:UniqueCover}
  Suppose that $\kappa$ is a composition of~$n$ such that
  $\vec\kappa=h(\lambda)$ and that $\t\in\Std(\lambda)_{w_\kappa}$ is a
  $w_\kappa$-transposable tableau. Then
  $\lambda_{\t,z}=\lambda^{(z)}$, for $1\le z\le d$. In particular,
  $X_{\t,z}$ is connected.
\end{Corollary}

\begin{proof}
  By \cref{L:covers} there is a unique sequence
  of self-conjugate partitions $\lambda^{(0)},\dots,\lambda^{(d)}$ such that
  $|\lambda^{(z)}|=\kappa_1+\dots+\kappa_z$ and
  $\lambda^{(z)}/\lambda^{(z-1)}$ is connected, for $1\le z\le d$. On the
  other hand, if $\t$ is a $w_\kappa$-transposable $\lambda$-tableau then
  $\lambda_{\t,1},\dots,\lambda_{\t,d}$ is a symmetric covering of
  $\lambda$. Therefore, $\lambda_{\t,z}=\lambda^{(z)}$,
  for $1\le z\le d$.
\end{proof}

By \cref{C:UniqueCover}, if $\t$ is a $w_\kappa$-transposable tableau then
the partitions $\lambda_{\t,z}$ depends only on~$\lambda$,  $\kappa$ and~$z$ and
not on~$\t$. Let $\t_{\kappa,z}$ be the restriction of~$\t$ to
$\lambda^{(z)}/\lambda^{(z-1)}$.
Then the entries of~$\t_{z,\kappa}$ are precisely the integers
$\set{k_z+i|0\le i\le 2m_z}$. Let
\[ \Std(\lambda)_{\kappa,z}=\set{\t_{\kappa,z}|\t\in\Std(\lambda)_{w_\kappa}}. \]
Then $\Std(\lambda)_{\kappa,z}$ is the set of skew tableau of shape
$\lambda^{(z)}/\lambda^{(z-1)}$ that are $w_\kappa$-transposable in the obvious
sense. The sets $\Std(\lambda)_{\kappa,z}$ are compatible with the different
$\sim_y$-equivalence classes in the sense that if $\s\sim_z\t$ then
$\s_{\kappa,y}=\t_{\kappa,y}$ for $y\ne z$.  Hence, by \cref{C:UniqueCover},
\begin{equation}\label{E:SimEquivClasses}
\Std(\lambda)_{w_\kappa}\bijection
        \prod_{z=1}^{d(\lambda)}\Std(\lambda)_{\kappa,z},
\end{equation}
where the bijection is given by $\t\mapsto(\t_{\kappa,1},\dots,\t_{\kappa,d})$.
Abusing notation slightly, we restrict the equivalence relation $\sim_z$ to
$\Std(\lambda)_{\kappa,z}$ and let $\Std[\lambda]_{\kappa,z}$ be the set of
$\sim_z$-equivalence classes of skew tableaux in $\Std(\lambda)_{\kappa,z}$.

Recall that we have fixed a labelling $\{x_0,\dots,x_{m_z}\}$ of~$X_\T$ such
that if $x_i\lessdot x_j$ then $j\in\{i\pm1\}$. As $X_\T$ is a connected strip,
there are exactly two such labellings (unless $m_z=0$, of course). We fix the
labelling of $X_\T$ so that $x_i$ is the unique node in $X_\T$ such that
$c(x_i)=i$, for $0\le i\le m_z$.  In particular, $x_0=\omega=(z,z)$.

The next identity, which is quite cute, is the last piece of the puzzle that we
need to compute $\chi^\lambda(T_{w_\kappa}\tau)$.

\begin{Lemma}\label{L:CuteIdentity}
  Suppose that $\vec\kappa=h(\lambda)$ and $1\le z\le d(\lambda)$. Then
  \[\sum_{\T\in\Std[\lambda]_{\kappa,z}}\varepsilon_\T q^{2c_\T(m_z)}
           \prod_{i=0}^{m_z-1}\frac1{[c_\T(i+1)-c_\T(i)]}
           =q^{-m_z}\prod_{i=1}^{m_z}\frac{[2i]}{[2i-1]}.\]
\end{Lemma}

\begin{proof} For clarity we write $m=m_z$ throughout the proof. We want to argue
  by induction on~$m$ but, before we can do this, we first need
  a more explicit description of $\Std[\lambda]_{\kappa,z}$. Recall
  the sign sequence $\varepsilon_\T$ from \cref{D:SimEquiv} and set
  $\underline\varepsilon_\T=(\varepsilon_\T(x_0),\dots,\varepsilon_\T(x_{m}))$.
  We claim that the map $\T\mapsto\underline\varepsilon_\T$ defines a bijection
  \[\Std[\lambda]_{\kappa,z}\bijection
  \Em=\set{(1,\varepsilon_1,\dots,\varepsilon_{m})|\varepsilon_i=\pm1}.\]
  (Consequently, $\#\Std[\lambda]_{\kappa,z}=2^m$.)
  To see this first note that if $\T$ is a $\sim$-equivalence class then
  $\varepsilon_\T(x_0)=1$. Therefore, by \cref{D:SimEquiv}(c) and
  \cref{C:UniqueCover} there are at most $2^{m}$ different~$\sim$-equivalence
  classes of tableaux in $\Std[\lambda]_{\kappa,z}$. On the other hand, if $\t$
  is any $w_\kappa$-transposable tableau then we obtain~$2^{m}$ different
  $w_\kappa$-transposable tableaux, each of which is in a different
  $\sim$-equivalence class, by swapping all pairs $\{\t(x),\t(x')\}$ of
  diagonally opposite entries in~$\t$ in all possible ways, for
  $x\in X_\T\setminus\{\omega\}$.  This establishes the claim.

  Identifying $\Std[\lambda]_{\kappa,z}$ and $\Em$ via the bijection
  above, we now argue by induction on~$m$.  If $m=0$ then both sides of
  the identity that we want to prove are equal to~$1$, so there is
  nothing to prove (by convention empty products are~$1$). We now assume
  that~$m>1$ and, by induction, that the lemma holds for~$\Em$. Let $\T$
  and $\T'$ be equivalence classes in $\Em[m+1]$ with
  $\varepsilon_{\T'}(x_r)=\varepsilon_{\T}(x_r)$, for $0\le r\le m$, and
  $\varepsilon_{\T}(x_{m+1})=1=-\varepsilon_{\T'}(x_{m+1})$. The
  contribution that $\T$ and $\T'$ make to the sum over~$\Em[m+1]$ is
  \[
     \varepsilon_\T q^{2c_\T(m)}\Biggl(
        \frac{q^{2(c_\T(m+1)-c_\T(m))}}{[c_\T(m+1)-c_\T(m)]}
        -\frac{q^{2(-c_\T(m+1)-c_\T(m))}}{[-c_\T(m+1)-c_\T(m)]}\Biggr)
        \prod_{i=0}^{m-1}\frac1{[c_\T(i+1)-c_{\T}(i)]}.
  \]
  By \cref{C:UniqueCover}, the posets $X_\T$ and $X_{\T'}$ are both
  connected so $c_\T(m+1)=m+1=-c_{\T'}(m+1)$ and
  $c_\T(m)=c_{\T'}(m)=\varepsilon_T(x_m)m=\pm m$.
  If $\varepsilon_{\T}(x_{m})=1$ then $c_\T(m)=m$ and
  \[ \frac{q^{2(m+1-c_\T(m))}}{[m+1-c_\T(m)]}
    -\frac{q^{-2(1+m+c_\T(m))}}{[-1-m-c_\T(m)]}
    =\frac{q^2}{[1]}-\frac{q^{-2(1+2m)}}{[-1-2m]}
    =q^{-1}\frac{[2m+2]}{[2m+1]}.
  \]
  Similarly if $\varepsilon_{\T}(x_{m})=-1$ then essentially the same
  calculation shows that the difference of these two terms is again
  $q^{-1}[2m+2]/[2m+1]$. All of the sequences in $\Em[m+1]$ occur in
  pairs of the form $(\underline\varepsilon_\T,\underline\varepsilon_{\T'})$, so
  \begin{align*}
    \sum_{\underline\varepsilon_\T\in\Em[m+1]}\varepsilon_\T q^{c_\T(m+1)}
      \prod_{i=0}^{m}\frac1{[c_\T(i+1)-c_\T(i)]}
      &=q^{-1}\frac{[2m+2]}{[2m+1]}\sum_{\SS\in\Em[m-1]}\varepsilon_{\SS} q^{c_{\SS}(m)}
      \prod_{i=0}^{m}\frac1{[c_{\SS}(i+1)-c_{\SS}(i)]}\\
      &=q^{-(m+1)}\prod_{i=1}^{m+1}\frac{[2i]}{[2i-1]}
  \end{align*}
  where the last equality follows by induction. This completes the proof
  of the inductive step and the lemma.
\end{proof}

We can now prove \cref{T:TauCharacters}, which is arguably the most
important result in this paper.

\begin{proof}[Proof of \cref{T:TauCharacters}]
  Recall that $h(\lambda)=(h_1,\dots,h_{d(\lambda)})$.
  We need to show that
  \[\chi^\lambda(T_{w_\kappa}\tau)
         =\begin{cases}\displaystyle
           \epsilon_\kappa(-\sqrt{-1})^{\half(n-d(\lambda))} q^{-\half n}
             \prod_{i=1}^{d(\lambda)}\sqrt{[h_{i}]},&
                 \text{if }\vec\kappa=h(\lambda),\\
           0,&\text{otherwise,}
       \end{cases}
  \]
  The results in this section show that
  $\chi^\lambda(T_{w_\kappa}\tau)=0$ if $\kappa$ is a composition
  of~$n$ such that $\vec\kappa\ne h(\lambda)$; see
  \cref{C:odd}, \cref{C:Many}, \cref{C:even} and \cref{C:zero!}.
  It remains to compute $\chi^\lambda(T_{w_\kappa}\tau)$ when
  $\vec\kappa=\lambda$.  Using \cref{C:transposable},
  \cref{E:SimEquivClasses} and \cref{P:transposableReduction},
  \begin{align*}
    \chi^\lambda_\H(T_{w_\kappa}\tau)
    &=\sum_{\t\in\Std(\lambda)_{\kappa,z}}\prod_{z=1}^{d(\lambda)}
           \gamma_{\t,z}(w_\kappa)
     =\prod_{z=1}^{d(\lambda)}\sum_{\T\in\Std[\lambda]_{\kappa,z}}\sum_{\t\in\T}
           \gamma_{\t,z}(w_\kappa)\\
    &=\prod_{z=1}^{d(\lambda)}\sum_{\T\in\Std[\lambda]_{\kappa,z}}
           (-1)^{m_z}q^{2c_\T(m_z)}\varepsilon_\T\epsilon_{X_\T}\alpha(X_\T)\,
             \prod_{i=0}^{m_z-1}\frac1{[c_\T(i+1)-c_\T(i)]}.
  \end{align*}
  Fix $z$ with $1\le z\le d(\lambda)$. By \cref{C:UniqueCover}, if
  $\T\in\Std[\lambda]_{\kappa,z}$ then~$X_\T=X_z$ depends
  only on $\kappa$ and $z$ (and $\lambda$). Similarly,
  the sign $\epsilon_{X_z}=\epsilon_{X_\T}$ depends only on~$\kappa$ and~$z$.
  Expanding the definition of $\alpha(X_\T)$, from \cref{E:alphaT},
  \[\alpha(X_\T)=\prod_{x\in X_\tau\setminus\set{\omega}}\alpha_{2c(x)}
               =\prod_{=0}^{m_z}\alpha_{2j}
               =(\sqrt{-1})^{m_z}\frac{\sqrt{[2m_z+1]}}{\sqrt{[1]}}
                \prod_{=1}^{m_z}\frac{[2i-1]}{[2i]}.
  \]
  Consequently, $\alpha(X_z)=\alpha(X_\T)$ also depends only on~$\kappa_z$.
  Set $\epsilon'_\kappa=\prod_{z=1}^{d(\lambda)}\epsilon_{X_z}$.
  Then the equation above becomes
  \begin{align*}
    \chi^\lambda_\H(T_{w_\kappa}\tau)
       &=\epsilon'_\kappa\prod_{z=1}^{d(\lambda)}(-1)^{m_z}\alpha(X_\T)
          \sum_{\T\in\Std[\lambda]_{\kappa,z}}
           \varepsilon_\T q^{2c_\T(m_z)}\,
             \prod_{i=0}^{m_z-1}\frac1{[c_\T(i+1)-c_\T(i)]}\\
       &=\epsilon'_\kappa\prod_{z=1}^{d(\lambda)}(-1)^{m_z}\alpha(X_\T)
            q^{-m_z}\prod_{i=1}^{m_z}\frac{[2i]}{[2i-1]},
            \hspace*{20mm}\text{by \cref{L:CuteIdentity}},\\
       &=\epsilon'_\kappa(-q^{-1}\sqrt{-1})^{m_1+\dots+m_z}
         \prod_{z=1}^{d(\lambda)}\frac{\sqrt{[2m_z+1]}}{\sqrt{[1]}},
  \end{align*}
  where the last equality uses the formula for $\alpha(X_z)$ given
  above. Now, $\sqrt{[1]}=\sqrt{q}$,
  $\sum_{z=1}^{d(\lambda)}m_z=\half(n-d(\lambda))$ and
  $\set{h_1,\dots,h_{d(\lambda)}}=\set{2m_z+1|1\le z\le d(\lambda)}$.
  Therefore, we have shown that
  \[ \chi^\lambda_\H(T_{w_\kappa}\tau)
          =\epsilon'_\kappa(-\sqrt{-1})^{\half(n-d(\lambda))}
              q^{-\half n}\prod_{i=1}^{d(\lambda)}\sqrt{[h_i]}.
  \]
  To complete the proof it remains to show that
  $\epsilon'_\kappa=\epsilon_\kappa$, where~$\epsilon_\kappa$ is the
  sign defined in \cref{T:TauCharacters}. By \cref{E:PosetSign},
  $\epsilon_{X_z}=(-1)^{\rho(z)}$ where $\rho(z)=\#\set{0\le i<m_z|x_{i+1}>x_i}$.
  Therefore, $\rho(z)+1$ is equal to the number of different
  rows in $\diag(\lambda)\cap X_z$. More precisely, if $x_{m_z}=(r,c)$ then
  $\rho(z)=z-r$. Armed with this observation it follows that
  that if~$\kappa_y<\kappa_{y+1}$ then $\epsilon'_\kappa=-\epsilon'_{\pi}$ where
  $\pi=(\kappa_1,\dots,\kappa_{y+1},\kappa_y,\dots,\kappa_d)$ is obtained from
  $\kappa$ by swapping $\kappa_y$ and $\kappa_{y+1}$. It follows that
  $\epsilon'_\kappa=(-1)^{\#\set{1\le y<z\le d(\lambda)|\kappa_y<\kappa_z}}
                   =\epsilon_\kappa$
  as required. The theorem is proved.
\end{proof}

\section{Computing all character values of $\An$}\label{S:characters}

  Combining \cref{C:NonSplitCharacters}, \cref{C:SplitCharacters} and
  \cref{T:TauCharacters} we can compute the values of all of the
  irreducible characters of~$\An$ on the elements $A_{w_\kappa}$,
  whenever $\kappa$ is a composition of~$n$.  Now
  $\set{w_\kappa|\kappa\in\Parts}$ is
  a set of minimal length conjugacy class representatives for
  $\Sym_n$, however, it is not a complete set of minimal length
  conjugacy class representatives for $\Alt_n$ because it does not
  contain the elements $w_\kappa^-$ whenever $\kappa$ is a partition of
  $n$ with distinct odd parts. Consequently, if $\lambda=\lambda'$ then
  we do not yet know the value of the characters
  $\chi^{\lambda\pm}(A_{w_\kappa^-})$ when $\kappa$ is a composition of
  $n$ with distinct odd parts. More importantly, if $\chi$ is a
  character of~$\An$ then we do not known how to compute $\chi(A_w)$ if
  $w\ne w^\pm_\kappa$ for some composition~$\kappa$. This section shows
  that the characters of~$\An$ are determined by characters of~$\Hn$ and
  \cref{T:TauCharacters}.

  The conjugacy classes of $\Sym_n$ and $\Alt_n$ have been described in
  \cref{S:conjugacy}. It is well-known that the characters of $\Hn$ are
  not class functions in the sense that if $\chi$ is a character then,
  in general, $\chi(T_v)$ and $\chi(T_w)$ are not necessarily equal if
  $v$ and $w$ are conjugate in $\Sym_n$. For example, this already
  happens for the character of the trivial character because
  $1_\H(T_w)=q^{\len(w)}$, for $w\in\Sym_n$. Nonetheless, as we now
  recall, the characters of~$\Hn$ are determined by their values on a
  set of minimal length conjugacy class representatives.

  For each conjugacy class $C\in\CC(\Sym_n)$ let
  $\Cmin=\set{x\in C|\ell(x)\le\ell(y)\text{ for all }y\in C}$ be the
  set of minimal length conjugacy class representatives. Fix an element
  $w_C\in\Cmin$ for each conjugacy class. For example, if~$C=C_\kappa$
  is the $\Sym_n$-conjugacy class of permutations of cycle type
  $\kappa\in\Parts$ then we could set $w_C=w_\kappa$.

  \begin{Theorem}[\protect{Geck-Pfeiffer \cite[\S8.2]{GeckPfeiffer:book}}]
    \label{T:GeckPfeiffer}
    Suppose that $\chi$ is a character of $\Hn$. Then there exist Laurent
    polynomials
    $\set{f_{C,w}(q)\in\Z[q,q^{-1}]|w\in\Sym_n\text{ and }C\in\CC(\Sym_n)}$,
    which do not depend,on~$\chi$, such that
    \[ \chi(T_w) = \sum_{C\in\CC(\Sym_n)} f_{C,w}(q)\chi(T_{w_C}) \]
    for any character $\chi$ of~$\Hn$. Moreover, if $C\in\CC(\Sym_n)$
    then $f_{C,w}(q)=0$ if $\ell(w)<\len(w_C)$ and $\chi(T_x)=\chi(T_{w_C})$
    for all $x\in\Cmin$.
  \end{Theorem}

  The polynomials $\{f_{C,w}(q)\}$ are the \textbf{class polynomials} of~$\Hn$.

  Recall the basis $\set{B_w|w\in\Alt_n}$ of $\An$ from \cref{P:BBasis}
  together with the involutions $\bar$, on~$\Zcal$ and $\epsilon$,
  on~$\H_{\Zcal,q}(\Sym_n)$, from \cref{C:KLInvariance}. Extend these involutions
  to~$\F$ and $\Hn=\H_{\F,q}(\Sym_n)$ by setting
  \[\overline{\sqrt{-1}}=\sqrt{-1}\quad\text{and}\quad
    \overline{\sqrt{[k]}}=\sqrt{-1}\sqrt{[-k]}, \qquad\text{for }k>0.
  \]
  The argument of \cref{C:KLInvariance} applies without change
  over~$\F$.

  \begin{Lemma}
    Suppose that $\lambda\in\Parts$ and $w\in\Sym_n$.
    Then
    $\chi^\lambda(B_w)=\varepsilon_w\overline{\chi^\lambda(B_w)}$.
  \end{Lemma}

  \begin{proof} Write $v_s B_w=\sum_\t b_{\s\t}(w)v_\s$, for some
    $b_{\s\t}(w)\in\Zcal$. Now $\epsilon(B_w)=\varepsilon_wB_w$, by
    \cref{C:KLInvariance}, so we must have
    $b_{\s\t}(w)=\varepsilon_w\overline{b_{\s\t}(w)}$. Taking
    traces proves the lemma.
      \end{proof}

  \begin{Lemma}\label{L:BVanishing}
    Suppose that $\lambda\in\Parts$ and $w\in\Sym_n\setminus\Alt_n$.
    Then $(\chi^\lambda+\chi^{\lambda'})(B_w)=0$.
  \end{Lemma}

  \begin{proof}
    By \cref{C:Hash}, $\chi^{\lambda'}(B_w)=\chi^\lambda(B_w^\#)$
    and $B_w^\#=-B_w$ since $w\notin\Alt_n$. Therefore,
    \[(\chi^\lambda+\chi^{\lambda'})(B_w)
          =\chi^\lambda(B_w)+\chi^{\lambda'}(B_w)
          =\chi^\lambda(B_w)+\chi^\lambda(B_w^\#)
          =\chi^\lambda(B_w)-\chi^\lambda(B_w)
          =0.
    \]
  \end{proof}

  The analogues of the last two lemmas are false for the $A$-basis
  of~$\Hn$.

  We can now prove an analogue of the Geck-Pfeiffer theorem for the
  irreducible characters of $\An$ that are indexed by partitions that
  are not self-conjugate.

  \begin{Proposition}\label{P:BGeckPfeiffer}
  Suppose that $\lambda\in\Parts$ and $w\in\Alt_n$. Then there exist
  polynomials $g_{C,w}(q)\in\Zcal$, which are independent of~$\lambda$, such
  that
  \[ \chi^\lambda_\A(A_w)
        =\sum_{C\in\CC(\Alt_n)}g_{C,w}(q)\chi^\lambda_\A(A_{w_C}).
  \]
  Moreover, $\chi^\lambda_\A(A_x)=\chi^\lambda_\A(A_{w_C})$ whenever
  $x\in\Cmin$.
  \end{Proposition}

  \begin{proof}
    We argue by induction on $\len(w)$. If $\ell(w)=1$ then $w=1$ and
    $\chi^\lambda_\A(A_w)=\dim S(\lambda)$ and there is nothing to prove.
    Suppose then that $\ell(w)>0$ and that $w\in D\in\CC(\Sym_n)$. Recall
    from \cref{C:NonSplitCharacters} that
    \[
        \chi^\lambda_\A(A_w)=\half\(\chi^\lambda(T_w)+\chi^{\lambda'}(T_w)\).
    \]
    Therefore, if $w\in\Cmin[D]$ then
    $\chi^\lambda_\A(A_w)=\chi^\lambda_\A(A_{w_{D}})$ by
    \cref{T:GeckPfeiffer}. In particular,
    $\chi^\lambda_\A(A_w)=\chi^\lambda_\A(A_{w_{D}})$ whenever
    $w\in\Cmin[D]$ (and, in fact, $\chi^\lambda_\A(w)$ depends only on the
    $\Sym_n$ conjugacy class of~$w$ and not the $\Alt_n$-conjugacy
    class). Hence, the  proposition follows when $w\in\Cmin[D]$
    since $D$ is a disjoint union of $\Alt_n$ conjugacy classes.

   Now suppose that $w$ is not of minimal length in its conjugacy class.
   Applying \cref{T:GeckPfeiffer},
   \begin{align*}
     \chi^\lambda_\A(A_w)
     &=\sum_{C\in\CC(\Sym_n)}\tfrac12{f_{C,w}}\(\chi^\lambda(T_{w_C})
                   +\chi^{\lambda'}(T_{w_C})\).\\
    \intertext{For $y\in\Sym_n$ we can write $T_v=\sum_{y\le v}s_{yv}B_y$, for some
    $s_{yv}\in\Zcal$. Therefore,}
     \chi^\lambda_\A(A_w)
         &=\sum_{y\in\Sym_n}\Big(\sum_{\substack{C\in\CC(\Sym_n)\\y\le w_C}}
         \tfrac12{s_{yw_C}f_{C,w}}\Big)(\chi^\lambda+\chi^{\lambda'})(B_y)\\
         &=\sum_{y\in\Alt_n}\Big(\sum_{\substack{C\in\CC(\Sym_n)\\y\le w_C}}
         s_{yw_C}f_{C,w}\Big) \chi_\A^\lambda(B_y)\\
     \intertext{where the last equality follows from \cref{L:BVanishing} and
     \cref{C:Hash}. Note that in the sum $y\in\Alt_n$.
     By \cref{P:BBasis}, if $y\in\Alt_n$ then $B_y=\sum_{x\le
     y}r_{xy}A_x$, for some $r_{xy}\in\Zcal$ such that $r_{xy}\ne0$ only
     if $x\in\Alt_n$. Hence, the sum above becomes}
     \chi^\lambda_\A(A_w)
       &=\sum_{x\in\Alt_n}\Bigg(\sum_{\substack{y\in\Alt_n\\x\le y}}
            \sum_{\substack{C\in\CC(\Sym_n)\\y\le w_C}}
               r_{xy}s_{yw_C}f_{C,w}\Bigg) \chi_\A^\lambda(A_x).
    \end{align*}
    Since $w$ is not of minimal length in its conjugacy class
    $f_{C,w}\ne0$ only if $\ell(w_C)<\ell(w)$. Consequently,
    $\chi^\lambda_\A(A_x)$ contributes to $\chi^\lambda_\A(A_w)$ only if
    $\ell(x)\le\ell(y)\le\ell(w_C)<\ell(w)$ in the last displayed
    equation. Therefore, by induction on length, $\chi^\lambda_\A(A_x)$
    can be written in the required form. It follows that
    $\chi^\lambda_\A(A_w)$ can be written in the required form, so
    the proof is complete.
  \end{proof}

  \cref{P:BGeckPfeiffer} implies that if $\chi$ is a character of~$\Hn$
  then the restriction of~$\chi$ to~$\An$ is determined by the character
  values $\set{\chi(A_{w_C})|C\in\CC(\Alt_n)}$ on the minimal length
  conjugacy classes of~$\Alt_n$. In particular, if~$\lambda\ne\lambda'$
  then $\chi^\lambda_\A$ is an irreducible character of~$\An$ and
  $\chi^\lambda_\A(A_w)$ is determined by \cref{P:BGeckPfeiffer}, for
  $w\in\Alt_n$.

  We now consider the characters of~$\An$ that are not the restriction
  of a character of~$\An$. For this it is enough to consider the
  irreducible characters $\chi^{\lambda\pm}_\A$ of~$\An$, where
  $\lambda$ is a self-conjugate partition of~$n$. To compute
  $\chi^{\lambda\pm}_\A(A_w)$, for $w\in\Alt_n$, we need to delve deeper
  into the proof of \cref{T:GeckPfeiffer}.

  Following \cite{GeckPfeiffer:irred}, let $\longrightarrow$ be the
  transitive closure of the relation $\xxrightarrow{s}$ on~$\Sym_n$, for $s\in S$
  where if $x,y\in\Sym_n$ then $x\xxrightarrow{s} y$ if $y=sxs$ and
  $\ell(y)\le\ell(x)$. Secondly, let $\approx$ be the equivalence
  relation on~$\Sym_n$ generated by the relations
  $\Xrightarrow{u}$, for $u\in\Sym_n$, where
  $x\Xrightarrow{u} y$ if $\ell(x)=\ell(y)$ and either $ux=yu$ and
  $\ell(ux)=\ell(u)+\ell(x)$, or $xu=uy$ and $\ell(xu)=\ell(x)+\ell(u)$.

  We state the next important theorem  only in the special case of the
  symmetric groups but, using case-by-case arguments, Geck and
  Pfeiffer~\cite{GeckPfeiffer:irred} proved this result for any finite
  Coxeter group.  He and Nie~\cite{HeNie:CoxeterConjugacy} have
  generalised this result to the extended affine Weyl groups  using an
  elegant case-free proof.

  \begin{Theorem}[Geck and Pfeiffer~\cite{GeckPfeiffer:irred}]\label{T:Conjugacy}
      Let $C$ be a conjugacy class of~$\Sym_n$.
      \begin{enumerate}
          \item If $x\in C$ then there exists an element $y\in\Cmin$
          such that $x\longrightarrow y$.
          \item If $x,y\in\Cmin$ then $x\approx y$.
      \end{enumerate}
  \end{Theorem}

  We will use \cref{T:TauCharacters} and \cref{T:Conjugacy} to determine
  the irreducible characters of~$\An$.  Given \cref{T:Conjugacy}, it is
  straightforward to prove \cref{T:GeckPfeiffer}. In contrast, it will
  take quite a bit of work to compute the character values
  $\chi^\lambda(T_w\tau)$, for $w\in\Sym_n$. The complication is that we
  are not able to make use of \cref{T:Conjugacy}(b) because when we try to
  apply it to compute $\chi(T_w\tau)$ then terms $\chi(T_v\tau)$ with
  $\ell(v)>\ell(w)$ can appear and this breaks any argument that uses
  induction on the length of~$w$. Fortunately, \cref{T:TauCharacters}
  allows~$\kappa$ to be a composition, rather than just a partition, so we
  can replace \cref{T:Conjugacy}(b) with the next result.

  \begin{Lemma}\label{L:Cmin}
      Let $w\in\Cmin$, where $C$ is an $\Sym_n$-conjugacy class.
      Then $w\longrightarrow w_\kappa$, for some composition~$\kappa$
      of~$n$.
  \end{Lemma}

  \begin{proof}
     As remarked above, the elements of $\Cmin$ are Coxeter elements in
     some standard parabolic subgroup. When $x$ is written as a product
     of disjoint cycles the different cycles commute, so it is enough to
     show if $x$ is a Coxeter element in $\Sym_{n+1}$ then
     $x\longrightarrow w_{(n+1)}=s_1\dots s_n$. As $x$ is a Coxeter
     element in~$\Sym_{n+1}$ we can write $x=s_{i_1}\dots s_{i_n}$,
     where $\set{i_1,\dots,i_n}=\set{1,\dots,n}$.  For $r=0,1,\dots n$
     we claim that
     $x\longrightarrow s_1\dots s_{r-1}s_{j_r}\dots s_{j_n}$ for some
     $j_r,\dots,j_n$ such that and $\set{j_r,\dots,j_n} =\set{r,r+1,\dots,n}$.
     To prove this we argue by induction on~$r$. If $r=0$ then
     there is nothing to prove because we can take $j_t=i_t$, for $1\le
     t\le n$. Hence, by induction, we may assume that $x\longrightarrow
     s_1\dots s_{r-1}s_{j_r}\dots s_{j_n}$ for suitable $j_r,\dots,j_n$.
     By assumption, $j_k=r$ for some $k\ge r$. If $k=r$ the inductive
     assumption automatically holds. If $k>r$ then, since
     $s_{j_r},\dots,s_{j_{k-1}}$ commute with $s_1,\dots,s_{r-1}$,
     \begin{align*}
       x &\xxrightarrow{s_{j_r}} s_{j_r}xs_{j_r}
          =s_1\dots s_{r-1} s_{j_{r+1}}\dots s_rs_{j_{k+1}}\dots s_{j_n}s_{j_r}\\
          &\xxrightarrow{s_{j_{r+1}}}s_{j_{r+1}}s_{j_r}xs_{j_r}s_{j_{r+1}}
            \xxrightarrow{s_{j_{r+2}}}\dots\\
          &\xxrightarrow{s{j_{k-1}}}
              s_{j_{k-1}}\dots s_{j_r}xs_{j_r}\dots s_{j_{k-1}}
          =s_1\dots s_r s_{j_{k+1}}\dots s_{j_n}s_{j_r}\dots s_{j_{k-1}},
     \end{align*}
     proving the inductive step. The proves the claim and hence the lemma.
  \end{proof}

  We need some preparatory lemmas before we can apply these results to
  compute $\chi^\lambda(A_w\tau)$, for an arbitrary permutation
  $w\in\Alt_n$. Recall that~$\tau$ is the unique $\F$-linear
  endomorphism of~$S(\lambda)$ such that $v_\t\tau=v_{\t'}$, for
  $\t\in\Std(\lambda)$.

  \begin{Lemma}\label{L:TauCommuting}
  Suppose that $\lambda$ is a self--conjugate partition and
  $w\in\Sym_n$. Then
  \[ v\tau T_w=vT_w^\#\tau\qquad\text{and}\qquad vT_w\tau=v\tau T_w^\#, \]
  for all $v\in S(\lambda)$.
  \end{Lemma}

  \begin{proof} By linearity, it is enough to consider the case when
    $v=v_\t$, for $\t\in\Std(\lambda)$. Then  $v_\t T_w^\#=v_\t\tau T_w\tau$,
    by \cref{twist}, implying the result.
  \end{proof}

  Recall that $S=\{s_1,\dots,s_{n-1}\}$ is the set of Coxeter generators of $\Sym_n$.
  If $s\in S$ and $w\in\Sym_n$ then $\len(sws)-\len(w)\in\{0,\pm2\}$. We need
  the following well-known result that is proved, for example, in
  \cite[Lemma~1.9]{M:Ulect}.

  \begin{Lemma}
    Suppose that $s\in S$ and $w\in\Sym_n$ such that $\len(sws)=\len(w)$ and
    $\len(sw)=\len(ws)$. Then $w=sws$.
  \end{Lemma}

  Hence, if $sws\ne w$ and $\ell(w)=\ell(sws)$ then either
  $\ell(ws)>\ell(w)=\ell(sws)>\ell(sw)$ or $\ell(sw)>\ell(w)=\ell(sws)>\ell(ws)$.
  The next lemma relates the character values
  $\chi^\lambda(T_x\tau)$ and $\chi^\lambda(T_y\tau)$ whenever $x\longrightarrow y$.
  It would be better if we could prove an analogue of this result for
  the characters of~$\An$, unfortunately, it is not clear how to do this.

  \begin{Lemma}\label{L:CharacterRelations}
  Suppose that $\lambda$ is a partition of~$n$ and suppose
  that $w\in\Sym_n$ and $s\in S$.
  \begin{enumerate}
  \item Suppose that $w\ne sws$ and $\ell(w)=\ell(sws)$. Then
  \[ \chi^\lambda(T_w\tau)=-\chi^\lambda(T_{sws}\tau)
                   +(q-q^{-1})\chi^\lambda(T_v\tau), \]
  where $v$ is uniquely determined by the requirements that $\ell(v)<\ell(w)$
  and $v\in\set{sw,ws}$.
  \item If $\ell(w)>\ell(sws)$ then
  $\chi^\lambda(T_w\tau) =-\chi^\lambda(T_{sws}\tau)$.
  \end{enumerate}
  \end{Lemma}

  \begin{proof}
  First consider~(a). Now $v=ws$  if and only if
  $\len(w)>\len(ws)$, in which case \cref{L:TauCommuting} implies that
  \begin{align*}
    \chi^\lambda(T_w\tau)&= \chi^\lambda(T_{ws}T_s\tau)
                          =\chi^\lambda(T_{ws}\tau T_s^\#)
    =\chi^\lambda(T_{ws}\tau(-T_s+q-q^{-1}))\\
    &=-\chi^\lambda(T_sT_{ws}\tau)+(q-q^{-1})\chi^\lambda(T_{ws}\tau),
  \end{align*}
  giving the result since $T_sT_{ws}=T_{sws}$.  The argument when
  $v=sw$ is similar.

  For part~(b), by \cref{L:TauCommuting},
  $\chi^\lambda(T_w\tau) = \chi^\lambda(T_sT_{sws}T_s\tau)
       = \chi^\lambda(T_{sws}\tau T_s^\#T_s) = -\chi^\lambda(T_{sws}\tau),$
  where the last equality follows because $T_s^\#=-T_s^{-1}$.
  \end{proof}

  Recall from \cref{S:conjugacy} that the $\Sym_n$-conjugacy class
  $C_\kappa$ splits into two $\Alt_n$-conjugacy classes if and only if
  the parts of~$\kappa$ are distinct and all odd. In this case $w^+_\kappa$
  and $w^-_\kappa=s_rw^+_\kappa s_r$, for some~$r$, are minimal length
  $\Alt_n$-conjugacy class representatives.

  \begin{Corollary}\label{C:KappaMinus}
      Suppose that $\lambda=\lambda'$ is a partition of~$n>1$ and that
      $\kappa$ is composition of~$n$ with $\vec\kappa= h(\lambda)$.
      Then
      \[\chi^\lambda(T_{w_\kappa^-})=-\chi^\lambda(T_{w_\kappa^+}).\]
  \end{Corollary}

  \begin{proof} By \cref{L:CharacterRelations}(a),
    $\chi^\lambda(T_{w^-_\kappa}\tau)=-\chi^\lambda(T_{w_\kappa}\tau)
                   +(q-q^{-1})\chi^\lambda(T_v\tau)$
    where $v$ is the unique permutation in $\set{s_rw^-_\kappa,w^-_\kappa s_r}$
    with $\ell(v)<\ell(w^-_\kappa)=\ell(w_\kappa)$. By
    \cref{T:Conjugacy}(a), if $C$ is the $\Sym_n$-conjugacy class
    containing~$v$ then $v\longrightarrow w_C$ for some $w_C\in\Cmin$.
    The~$\longrightarrow$ relation is generated by the two situations
    considered in \cref{L:CharacterRelations}, so it follows that
    $\chi^\lambda(T_v\tau)$ can be written as a $\Zcal$-linear
    combination of character values $\chi^\lambda(T_{w_D}\tau)$, where
    the sum is over $D\in\CC(\Sym_n)$ with $w_D\in\Cmin[D]$ and
    $\ell(w_D)\le\ell(v)$. By \cref{L:Cmin}, and possibly further
    applications of \cref{L:CharacterRelations}, we can assume that
    $w_D=w_\sigma$, for some composition~$\sigma$. If $\sigma$ is a
    composition appearing in this way then
    $\ell(w_\sigma)\le\ell(v)<\ell(w_\kappa)$, so $\vec\sigma\ne
    h(\lambda)$.  Therefore, $\chi^\lambda(T_{w_\sigma}\tau)=0$, by
    \cref{T:TauCharacters}, so that
    $\chi^\lambda(T_{w^-_\kappa}\tau)=-\chi^\lambda(T_{w_\kappa}\tau)$
    as required.
  \end{proof}

  More generally, if $w\in\Sym_n$ is a minimal length element of cycle
  type $h(\lambda)$ then
  $\chi^\lambda(T_w\tau)=\pm\chi^\lambda(T_{w_{h(\lambda)}}\tau)$.

  Using the argument of \cref{C:KappaMinus} we can now compute
  $\chi^\lambda(T_w\tau)$, for any $w\in\Sym_n$. However, if $\nu$ is a
  composition of~$n$ then this does not
  imply that $\chi^\lambda(T_{w^-_\nu}\tau)=\pm\chi^\lambda(T_{w^+_\nu}\tau)$
  because if $\ell(w_\nu)>\ell(w_{h(\lambda)})$ then applications of
  \cref{L:CharacterRelations} can introduce terms
  $\chi^\lambda(T_{w_\sigma}\tau)$ where~$\vec\sigma=h(\lambda)$. See
  the examples at the end of this section.

  We can now prove a stronger version of \cref{T:TauCharacters}. Let
  $\rightsquigarrow$ be the transitive relation on~$\Sym_n$ generated
  by~$w\rightsquigarrow v$ if either $w\longrightarrow v$
  or there exists $s\in S$ such that $\ell(w)=\ell(sws)$,
  $v\in\set{sw,ws}$ and $\ell(v)=\ell(w)-1$.

  \begin{Theorem}\label{T:TauCharactersII}
    Suppose that $\lambda=\lambda'$ and $w\in\Sym_n$. Let
    $h(\lambda)=(h_1,\dots,h_{d(\lambda)})$. Then
    there exists a polynomial $a^\lambda_w(a)\in\Z[(q-q^{-1})]$ such that
    $\deg_q a^\lambda_w(q)\le\ell(w)-\ell(w_{h(\lambda)})$ and
    \[\chi^\lambda(T_w\tau)=(-\sqrt{-1})^{\half(n-d(\lambda))}a^\lambda_w(q)q^{-\half n}
                    \prod_{i=1}^{d(\lambda)}\sqrt{[h_{i}]}.
    \]
    Moreover, $a^\lambda_w(q)\ne0$ only if there exists a composition $\kappa$
    such that $w\rightsquigarrow w_\kappa$ and~$\vec\kappa=h(\lambda)$.
  \end{Theorem}

  \begin{proof}
    As in the proof of \cref{C:odd}, by \cref{T:Conjugacy}(a) and
    repeated applications of \cref{L:CharacterRelations} and \cref{L:Cmin}, there
    exist polynomials $a_{\sigma w}(q)\in\Z[(q-q^{-1})]$ such that
    \[\chi^\lambda(T_w\tau)
        = \sum_\sigma a_{\sigma w}(q)\chi^\lambda(T_{w_\sigma}\tau),\]
    where the sum is over compositions $\sigma$ of~$n$ such that
    $w\rightsquigarrow w_\sigma$. Observe that the different cases of
    \cref{L:CharacterRelations} multiply the character values by
    either~$-1$ or $q-q^{-1}$, with $q-q^{-1}$ appearing only in the
    case of \cref{L:CharacterRelations}(a) where it arises as the
    coefficient of $\chi^\lambda(T_v\tau)$, where $v\le w$ and
    $\ell(v)=\ell(w)-1$. Therefore,
    $\deg_q a_{\sigma w}(q)\le\ell(w)-\ell(w_\sigma)$. Set
    $a^\lambda_w(q)=\sum_\kappa \epsilon_\kappa a_{\kappa w}(q)$, where
    the sum is over those compositions $\kappa$ with $\vec\kappa=h(\lambda)$. Then
    \[\chi^\lambda(T_w\tau) = a^\lambda_w(q)
               (-\sqrt{-1})^{\half(n-d(\lambda))}q^{-\half n}
                    \prod_{i=1}^{d(\lambda)}\sqrt{[h_{i}]}
    \]
    by \cref{T:TauCharacters}. All of the claims in the theorem now follow.
  \end{proof}

  The polynomials $a^\lambda_w(q)$ appearing in \cref{T:TauCharactersII}
  are, in principle, easy to compute.  Examples suggest that if~ $w\in
  C_\mu$ then $a^\lambda_w(q)\ne0$ only if $\mu\gedom\lambda$, where
  $\gedom$ is the dominance order on~$\Parts$.

  Combining \cref{P:SplitCharacters} and \cref{T:TauCharactersII} yields
  a weak analogue of the Geck-Pfeiffer \cref{T:GeckPfeiffer} for the
  irreducible characters of~$\An$ indexed by self-conjugate partitions.

\begin{Corollary}\label{T:BGeckPfeiffer}
  Suppose that $\lambda=\lambda'$ is a self-conjugate partition of~$n$
  and $w\in\Alt_w$. Then
  \[ \chi^{\lambda\pm}_\A(A_w)
        = \half\sum_{C\in\CC(\Alt_n)}g_{C,w}(q)\chi^\lambda_\A(A_{w_C})
        \pm\half (-\sqrt{-1})^{\half(n-d(\lambda))}a^\lambda_w(q)q^{-\half n}
            \prod_{i=1}^{d(\lambda)}\sqrt{[h_{i}]},
  \]
  where $g_{C,w}(q)$ and $a^\lambda_w(q)$ are the polynomials from \cref{P:BGeckPfeiffer}
  and \cref{T:TauCharactersII}, respectively.
\end{Corollary}

\begin{proof}
  By \cref{P:SplitCharacters} and \cref{T:TauCharactersII},
  \[ \chi^{\lambda\pm}_\A(A_w)
        = \half\chi^\lambda(T_w)
        \pm\half (-\sqrt{-1})^{\half(n-d(\lambda))}a^\lambda_w(q)q^{-\half n}
            \prod_{i=1}^{d(\lambda)}\sqrt{[h_{i}]}.
  \]
  On the other hand, using \cref{C:NonSplitCharacters},
  \[
   \chi^\lambda(T_w) =\half\bigl(\chi^\lambda(T_w)+\chi^{\lambda'}(T_w)\bigr)
                     =\chi^\lambda(A_w)
                     =\sum_{C\in\CC(\Alt_n)}g_{C,w}(q)\chi^\lambda_\A(A_{w_C}),
  \]
  where the last equality comes from \cref{P:BGeckPfeiffer}. The result
  follows.
\end{proof}

Therefore, if $\chi$ is a character of $\An$ then the values $\chi(A_w)$,
for $w\in\Alt_n$, can be written as a $\Zcal$-linear combination of the
character values $\chi(A_{w_C})$, where the sum is over the conjugacy
classes of~$\Alt_n$ and $w_C\in\Cmin$.

\cref{T:TauCharactersII} also allows us to improve on \cref{C:SplittingField}.

\begin{Corollary}
  Suppose that $\F$ is a field  of characteristic different from~$2$ that
  contains square roots
  \[ \sqrt{-1} \qquad\text{and}\qquad
     \sqrt{[2m+1]} \quad\text{ for }\quad 0\le m\le\tfrac{n-1}2.
  \]
  Then alternating Hecke algebra $\An[\F,q]$ is split semisimple.
\end{Corollary}

\begin{proof}
  By \cref{P:BGeckPfeiffer} and \cref{T:BGeckPfeiffer}, if~$\chi$
  is an irreducible character of~$\An$ then $\chi(a)\in\F$, for all $a\in\An$.
  As all of the characters of~$\An$ take values in~$\F$ it
  follows by general nonsense that~$\F$ is a splitting field for $\An$. See, for
  example,~\cite[7.15]{CurtisReiner:VolI}.
\end{proof}

  Finally, by way of example, we show that the obvious generalisations
  of \cref{T:GeckPfeiffer} to the characters of~$\An$ fail for both the
  $A$-basis and the $B$-basis of $\An$. This suggests that
  \cref{T:BGeckPfeiffer} (and \cref{P:BGeckPfeiffer}), may be the best
  results possible.

  \begin{Example}\label{Ex:FirstExample}
    We show that the characters values $\chi(A_w)$ are \textit{not}
    constant on the minimal length conjugacy class representatives of~$\Sym_n$. Take
    $\lambda=(3^3)$ so that $h(\lambda)=(5,3,1)$. Let
    $w=w_{(9)}=s_1s_2s_3s_4s_5s_6s_7s_8$. Then
    $\chi_\A^{\lambda\pm}(A_w)=\half\chi^\lambda(T_w)$
    by \cref{C:SplitCharacters} and \cref{T:TauCharacters}. Now consider
    $v=s_8s_5s_1s_2s_3s_4s_6s_7$. Then $v$ and $w$ are conjugate in~$\Sym_n$
    (although they are not conjugate in $\Alt_n$), and
    they have the same length. By \cref{C:SplitCharacters},
    \[\chi_\A^{\lambda\pm}(A_v)
       =\half\(\chi^\lambda(T_v)\pm\chi^\lambda(T_v\tau)\)
       =\half\(\chi^\lambda(T_w)\pm\chi^\lambda(T_v\tau)\),
    \]
    where the last equality comes from \cref{T:GeckPfeiffer} since $v$
    and $w$ are of minimal length in their conjugacy class. On the other
    hand, using the argument from the proof of \cref{T:TauCharactersII}
    shows that
    \[
      \chi^\lambda(T_v\tau)=(q-q^{-1})^2\chi^\lambda(T_1T_2T_3T_4T_6T_7\tau)
       =-q^{-4}(q-q^{-1})^2\sqrt{-1}\sqrt{[3]}\sqrt{[5]}.
    \]
    Consequently, $\chi_\A^{\lambda\pm}(A_v)\ne\chi_\A^{\lambda\pm}(A_w)$.
  \end{Example}

  \begin{Example}\label{Ex:SecondExample}
    Maintaining the notation of the last example with $\lambda=(3,3,3)$,
    brute force calculations using code written by the first author
    in \textsc{Sage}~\cite{sage} reveal the
    following:

    \[
      \begin{array}{ccc}\toprule
           x\in\Alt_9 & \chi^{\lambda\pm}_\A(A_x)
                      & \chi^{\lambda\pm}_\A(B_x)
            \\[2pt]\midrule
            w=s_1s_2s_3s_4s_5s_6s_7s_8 & 0 &2^8(q-q^{-1})^8\\
            v=s_8s_5s_1s_2s_3s_4s_6s_7
                  & \mp q^{-4}\sqrt{-1}\sqrt{[3]}\sqrt{[5]}(q-q^{-1})^2
                  & 2^8(q-q^{-1})^8\mp q^{-4}\sqrt{-1}\sqrt{[3]}\sqrt{[5]}(q-q^{-1})^2\\
       u=s_7s_8s_5s_1s_2s_3s_4s_6
                  & \pm q^{-4}\sqrt{-1}\sqrt{[3]}\sqrt{[5]}(q-q^{-1})^2
                  & 2^8(q-q^{-1})^8\pm q^{-4}\sqrt{-1}\sqrt{[3]}\sqrt{[5]}(q-q^{-1})^2
            \\\bottomrule
      \end{array}
    \]
    \medskip

    The three permutations $u$, $v$ and $w$ are of minimal length in
    their conjugacy class. All three elements are conjugate in~$\Sym_9$
    whereas only $u$ and $w$ are conjugate in~$\Alt_9$.  In particular,
    this calculation shows that the characters of~$\An$ are not
    constant on $A$ or $B$ basis elements indexed by minimal length conjugacy class
    representatives. Hence, the obvious generalisation of the
    Geck-Pfeiffer \cref{T:GeckPfeiffer} for the
    irreducible characters of~$\An$ is not true with respect to either the $A$ or
    $B$~bases of~$\An$.
  \end{Example}

\section*{Acknowledgements}
  This research was supported, in part, by the Australian Research Council. Some
  of the results in this paper appear in the second author's Ph.D.
  thesis~\cite{Ratliff:PhDThesis}.

%%%%%%%%%%%%%%%%%%%%%%%%%%%%%%%%%%%%%%%%%%%%%%%%%%%%%%%%%%%%%%%%
%\bibliography{papers}

\end{document}